\documentclass{article}
\usepackage[utf8]{inputenc} %
\usepackage[T1]{fontenc}    %
\usepackage{hyperref}       %
\usepackage{url}            %
\usepackage{booktabs}       %
\usepackage{amsfonts}       %
\usepackage{nicefrac}       %
\usepackage{microtype}      %
\usepackage{natbib}
\usepackage{fullpage}

\usepackage{color}

\def\x{{\bf{x}}}

\def\tilde{\hat}

\usepackage[utf8]{inputenc} %
\usepackage[T1]{fontenc}  %
\usepackage{etoolbox}
\newtoggle{icml}
\togglefalse{icml}
\newtoggle{aistats}
\togglefalse{aistats}
\usepackage{microtype}    %
\usepackage{comment}
\usepackage{abstract}     %
\usepackage{amsmath,amssymb,mathtools,amsfonts}
\usepackage[cal=boondoxo,scr=esstix]{mathalfa} %
\usepackage{lipsum}
\usepackage{amsthm}
\usepackage{nicefrac}       %
\usepackage{graphicx,color,wrapfig,epsfig,epstopdf}
\usepackage[font=small]{caption, subcaption}
\usepackage{enumitem}
\usepackage{booktabs}

\usepackage{algorithm,algorithmicx,algpseudocode}

\usepackage{hyperref}       %
\usepackage{cleveref}
\usepackage{url}            %
\usepackage{booktabs}       %
\usepackage{listings}
\usepackage{multirow}
\usepackage{stackengine}
\newcommand\numberthis{\addtocounter{equation}{1}\tag{\theequation}}

\allowdisplaybreaks[4]

\usepackage{thmtools}

\declaretheorem[thmbox=M,name=Theorem]{theorem}
\declaretheorem[name=Assumption]{assumption}
\declaretheorem[name=Remark]{remark}

\declaretheorem[thmbox=M,name=Lemma,sibling=theorem]{lemma}
\declaretheorem[thmbox=M,name=Corollary,sibling=theorem]{corollary}
\algblockdefx[BlockUntil]{BlockUntil}{EndBlockUntil}[1][]{\textbf{wait until}
  #1 \textbf{then}}{\textbf{end wait until}}

\def\cdummy{\cdot}
\def\backassign{=:}
\usepackage[T1]{fontenc}
\usepackage{microtype}
\usepackage{abstract}

\usepackage{etoolbox}
\usepackage[cal=boondoxo,scr=esstix]{mathalfa} %
\newenvironment{xiangru}
{\begingroup\color{black}}
{\endgroup\ignorespacesafterend}
\def\xr{\color{black}}
\setitemize{noitemsep,topsep=2pt,parsep=2pt,partopsep=0pt,leftmargin=0.02cm,itemindent=.5cm}
\setenumerate{noitemsep,topsep=2pt,parsep=2pt,partopsep=0pt,leftmargin=0cm,itemindent=.5cm}
\setdescription{noitemsep,topsep=2pt,parsep=2pt,partopsep=0pt,leftmargin=0.3cm,itemindent=0cm}

\usepackage{titlesec}
\titlespacing\section{0pt}{11pt plus 2pt minus 2pt}{5pt plus 1pt minus 2pt}
\titlespacing\subsection{0pt}{2pt plus 2pt minus 2pt}{0pt plus 1pt minus 2pt}
\titlespacing\subsubsection{0pt}{2pt plus 2pt minus 2pt}{0pt plus 1pt minus 2pt}
\titlespacing\paragraph{0pt}{0pt plus 2pt minus 2pt}{5pt plus 2pt minus 2pt}

\title{Revisit Batch Normalization: New Understanding from an Optimization View and a Refinement via
  Composition Optimization}

\usepackage{authblk}
 \author{Xiangru Lian}
 \author{Ji Liu}
 \affil[]{\texttt{admin@mail.xrlian.com, ji.liu.uwisc@gmail.com}}
 \affil[]{University of Rochester}

\begin{document}

\maketitle

\begin{abstract}
Batch Normalization (BN) has been used extensively in deep learning to achieve faster training process and better resulting models. However, whether BN works strongly depends on how the batches are constructed during training and it may not converge to a desired solution if the statistics on a batch are not close to the statistics over the whole dataset. In this paper, we try to understand BN from an optimization perspective by formulating the optimization problem which motivates BN. We show when BN works and when BN does not work by analyzing the optimization problem. We then propose a refinement of BN based on compositional optimization techniques called Full Normalization (FN) to alleviate the issues of BN when the batches are not constructed ideally. We provide convergence analysis for FN and empirically study its effectiveness to refine BN.
 \end{abstract}

\newlength{\firstfigurelength}
\newlength{\firstcaptionlength}
\setlength{\firstfigurelength}{0.45\textwidth}
\setlength{\firstcaptionlength}{0.5\textwidth}
\section{Introduction}

Batch Normalization (BN) \citep{ioffe2015batch} has been used extensively in
deep learning \citep{szegedy2016rethinking,he2016deep,silver2017mastering,huang2017densely,hubara2017quantized} to
accelerate the training process. During the training process, a BN layer
normalizes its input by the mean and variance computed within a mini-batch. Many
state-of-the-art deep learning models are based on BN such as ResNet
\citep{he2016deep,xie2017aggregated} and Inception \citep{szegedy2016rethinking,szegedy2017inception}. It is often
believed that BN can mitigate the exploding/vanishing gradients problem
  \citep{cooijmans2016recurrent} or reduce internal variance
  \citep{ioffe2015batch}. Therefore, BN has become a standard tool that is
implemented almost in all deep learning solvers such as Tensorflow
\citep{tensorflow2015-whitepaper}, MXnet \citep{chen2015mxnet}, Pytorch
\citep{paszke2017automatic}, etc.

Despite the great success of BN in quite a few scenarios, people with rich experiences with BN may also notice some issues with BN. Some examples are provided below. %

\begin{wrapfigure}{r}{\firstcaptionlength}
    \centering
    {\includegraphics[width=\firstfigurelength]{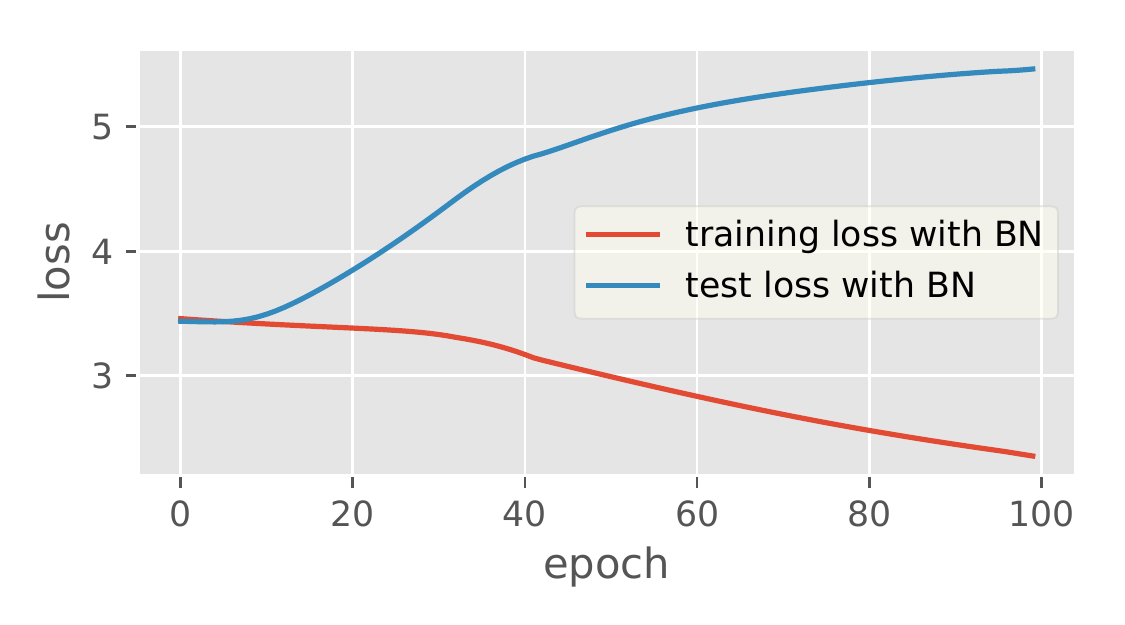}}
    \caption{\label{fig:toya} {(\bf Fail when batch size is 1)}
    Given a dataset comprised of three samples: ${[}0,0,0{]}$ with label $0$,
    ${[}1,1,1{]}$ with label $1$ and ${[}2,2,2{]}$ with label $2$, use the following
    simple network including one batch normalization layer, where the numbers in
    the parenthesis are the dimensions of input and output of a layer:
      linear layer (3 $\rightarrow$3) $\Rightarrow$ batch
      normalization %
      $\Rightarrow$ relu %
      $\Rightarrow$ linear layer (3 $\rightarrow$ 3) $\Rightarrow$ nll loss.
    \newline
    Train with batch size 1, and test on the same dataset. The test loss increases while the training loss
    decreases. %
    }
  \end{wrapfigure}

\paragraph{BN fails/overfits when the mini-batch size is 1 as shown in \Cref{fig:toya}}
We construct a simple network with 3 layers for a classification task. It fails
to learn a reasonable classifier on a dataset with only 3 samples as seen in \Cref{fig:toya}.

\begin{wrapfigure}{l}{\firstcaptionlength}
    \centering
    {\includegraphics[width=\firstfigurelength]{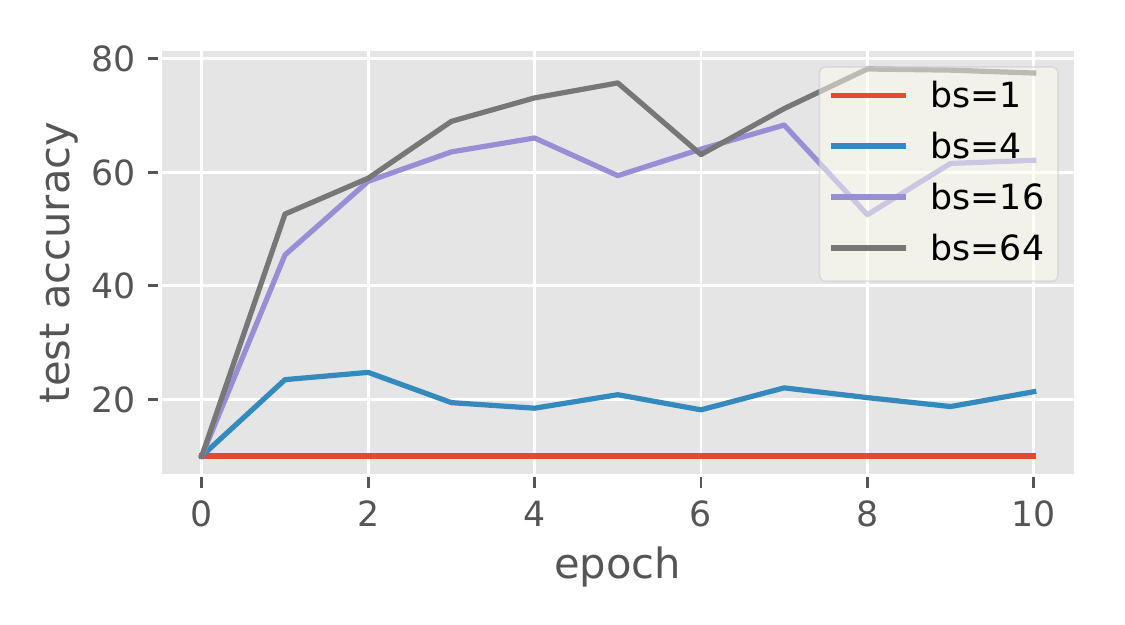}}
    \caption{\label{fig:cifar-bs} {(\bf Sensitive to the size of
        mini-batch)} The test accuracy for ResNet18 on CIFAR10
      dataset trained for 10 epochs with different batch sizes. The smaller the batch size
      is, with BN layers in ResNet, the worse the convergence result
      is.}
  \end{wrapfigure}

\begin{wrapfigure}{r}{\firstcaptionlength}
    \centering
    {\includegraphics[width=\firstfigurelength]{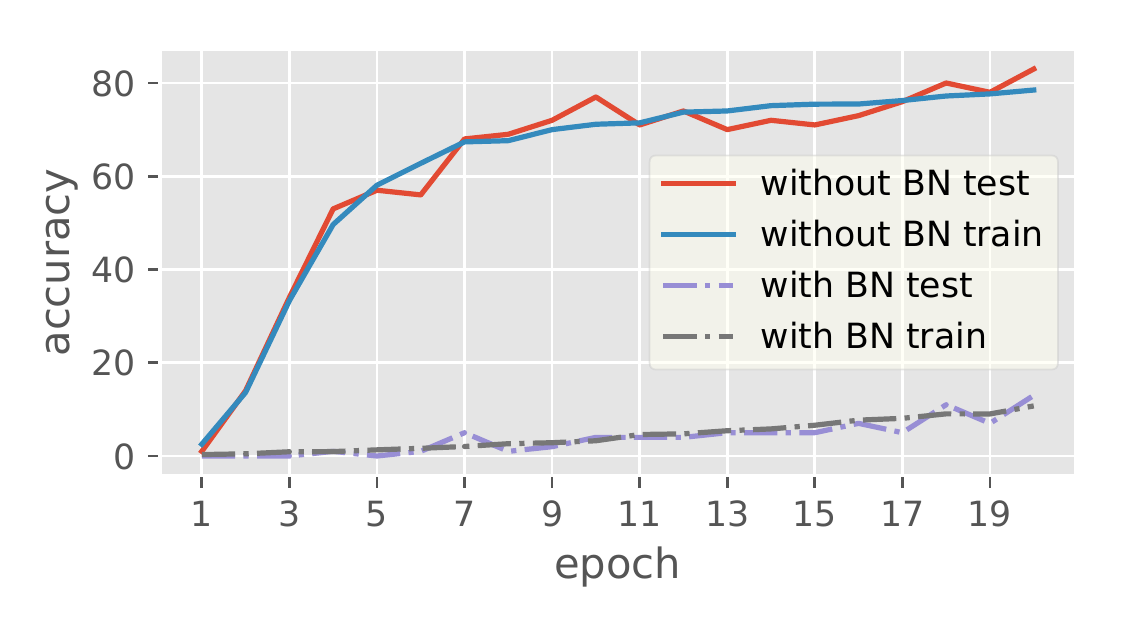}}
    \caption{
      \label{fig:logistic-bn} {(\bf Fail if data are with large variation)} The
      {\xr training and test} accuracy for a simple logistic regression problem
      on synthesized dataset with mini-batch size 20. We synthesize 10,000
      samples where each sample is a vector of 1000 elements. There are 1000
      classes. Each sample is generated from a zero vector by firstly randomly
      assigning a class to the sample (for example it is the $i$-th class), and
      then setting the $i$-th element of the vector to be a random number from 0
      to 50, and finally adding noise (generated from a standard normal
      distribution) to each element. We generate 100 test samples in the same
      way. A logistic regression classifier should be able to classify this
      dataset easily. However, if we add a BN layer to the classifier, the model
      no longer converges. }

  \end{wrapfigure}

  \paragraph{BN's solution is sensitive to the mini-batch size as shown in
    \Cref{fig:cifar-bs}} The test conducted in \Cref{fig:cifar-bs} uses ResNet18
  on the Cifar10 dataset. When the batch size changes, the neural network with
  BN training tends to converge to a different solution. Indeed, this
  observation is true even for solving convex optimization problems. It can also
  be observed that the smaller the batch size, the worse the performance of the
  solution.%
  \paragraph{BN fails if data are with large variation as shown in
    \Cref{fig:logistic-bn}} BN breaks convergence on simple convex logistic
  regression problems if the variance of the dataset is large.
  \Cref{fig:logistic-bn} shows the first 20 epochs of such training on a
  synthesized dataset. This phenomenon also often appears when using distributed
  training algorithms where each worker only has its local dataset, and the
  local datasets are very different.

Therefore, these observations may remind people to ask some fundermental questions:
\begin{enumerate}
\item Does BN always converge and/or where does it converge to?
\item Using BN to train the model, why does sometimes severe over-fitting
  happen?
\item Why can BN often accelerate the training process?
\item Is BN a trustable ``optimization'' algorithm?
\end{enumerate}

In this paper, we aim at understanding the BN algorithm from a rigorous
optimization perspective to show
\begin{enumerate}
\item BN always converges, but not solving either the objective in our common
  sense, or the optimization originally motivated.
\item The result of BN heavily depends on how to construct batches, and it can
  overfit predefined batches. BN treats the training and the inference
  differently, which makes the situation worse.
\item BN is not always trustable, especially when the batches are not
  constructed with randomly selected samples or the batch size is small.
\end{enumerate}

Besides these, we also provide the explicit form of the original objective BN
aims to optimize (but does not in fact), and propose a Multilayer Compositional
Stochastic Gradient Descent (MCSGD)
algorithm based on the compositional optimization technology to solve the
original objective. We prove the convergence of the proposed MCSGD algorithm and
empirically study its effectiveness to refine the state-of-the-art BN algorithm.

\section{Related work}

We first review {\em traditional normalization} techniques.
\cite{lecun1998efficient} showed that normalizing the input dataset makes
training faster. In deep neural networks, normalization was used before the
invention of BN. For example Local Response Normalization (LRU)
\citep{lyu2008nonlinear,jarrettk2009whatisthebestmulti,krizhevsky2012imagenet}
which computes the statistics for the neighborhoods around each pixel.

We then review {\em batch normalization} techniques. \cite{ioffe2015batch}
proposed the Batch Normalization (BN) algorithm which performs normalization
along the batch dimension. It is more global than LRU and can be done in the
middle of a neural network. Since during inference there is no ``batch'', BN
uses the running average of the statistics during training to perform inference,
which introduces unwanted bias between training and testing.
\cite{ioffe2017batch} proposed Batch Renormalization (BR) that introduces two
extra parameters to reduce the drift of the estimation of mean and variance.
However, Both BN and BR are based on the assumption that the statistics on a
mini-batch approximates the statistics on the whole dataset.

Next we review the {\em instance based normalization} which normalizes the input
sample-wise, instead of using the batch information. This includes Layer
Normalization \citep{ba2016layer}, Instance Normalization
\citep{ulyanov2016instance}, Weight Normalization \citep{salimans2016weight} and
a recently proposed Group Normalization \citep{wu2018group}. They all normalize
on a single sample basis and are less global than BN. They avoid doing
statistics on a batch, which work well when it comes to some vision tasks where
the outputs of layers can be of handreds of channels, and doing statistics on
the outputs of a single sample is enough. However we still prefer BN when it
comes to the case that the single sample statistics is not enough and more
global statistics is needed to increase accuracy.

Finally we review the {\em compositional optimization} which our refinement of
BN is based on. \cite{wang2016stochastic} proposed the compositional
optimization algorithm which optimizes the nested expectation problem $\min
\mathbb{E} f(\mathbb{E} g(\cdot))$, and later the convergence rate was improved
in \cite{wang2016accelerating}. The convergence of compositional optimization on
nonsmooth regularized problems was shown in \cite{huo2017accelerated} and a
variance reduced variant solving strongly convex
problems was analyzed in \cite{lian2016finite}. %
\section{Review the BN algorithm}

In this section we review the Batch Normalization (BN) algorithm \citep{ioffe2015batch}.

{\xr BN is usually implemented as an additional layer in neural
  networks. In each iteration of the training process, the BN layer normalizes
  its input using the mean and variance of each channel of the input batch to make its output having
  zero mean and unit variance. The mean and variance on the input batch is
  expected to be similar to the mean and variance over the whole dataset. In the inference process, the layer normalizes
  its input's each channel using the saved mean and variance, which are the running averages of
  mean and variance calculated during training.} This is described in
\Cref{algo:bn-train}\footnote{A linear transformation is often added after
  applying BN to compensate the normalization.} and \Cref{algo:bn-infer}.
  With BN, the input at the BN layer is normalized so that the next layer in the
  network accepts inputs that are easier to train on. In practice it has been
  observed in many applications that the speed of training is improved by using
  BN.

\begin{algorithm}
  \caption{Batch Normalization Layer (Training)\label{algo:bn-train}}
  \begin{algorithmic}[1]
    \Require Input $\mathbf{B}_{\text{in}}$, which is a batch of input.
    Estimated mean $\mu$ and variance $\nu$, and averaging constant $\alpha$.
    \State
    \begin{align*}
      \mu\gets& (1-\alpha)\cdot \mu + \alpha \cdot \textbf{\textit{mean}}(\mathbf{B}_{\text{in}}),\\
      \nu\gets& (1-\alpha)\cdot \nu + \alpha \cdot \textbf{\textit{var}}(\mathbf{B}_{\text{in}}).
    \end{align*}
    \Comment{$\textbf{\textit{mean}}(\mathbf{B}_{\text{in}})$ and
      $\textbf{\textit{var}}(\mathbf{B}_{\text{in}})$ calculate the mean and variance of
      $\mathbf{B}_{\text{in}}$ respectively. }

    \State {\bf Output} \[\mathbf{B}_{\text{out}} \gets \frac{\mathbf{B}_{\text{in}} -
      \textbf{\textit{mean}}(\mathbf{B}_{\text{in}})}{\sqrt{\textbf{\textit{var}}(\mathbf{B}_{\text{in}})
        +\epsilon}}.\]

    \Comment{$\epsilon$ is a small constant for numerical stability.}
  \end{algorithmic}
\end{algorithm}

\begin{algorithm}
  \caption{Batch Normalization Layer (Inference)\label{algo:bn-infer}}
  \begin{algorithmic}[1]
    \Require Input $\mathbf{B}_{\text{in}}$, estimated mean $\mu$ and variance
    $\nu$. A small constant $\epsilon$ for numerical stability.
    \State {\bf Output} \[\mathbf{B}_{\text{out}} \gets \frac{\mathbf{B}_{\text{in}} -
      \mu}{\sqrt{\nu +\epsilon}}.\]
  \end{algorithmic}
\end{algorithm}

\section{Training objective of BN}

In the original paper proposing BN algorithm \citep{ioffe2015batch},
authors do not provide the explicit optimization objective BN targets to solve. Therefore, many people may naturally think that BN is an optimization trick, that accelerates the training process but still solves the original objective in our common sense.
Unfortunately, this is not true! The actual objective BN solves is different from the objective in our common sense and also nontrivial.

\paragraph{Rigorous mathematical description of BN layer}
To define the objective of BN in a precise way, we need to define the normalization operator
$f^{B, \sigma}_W$ that maps a \emph{function} to a \emph{function} associating
with a mini-batches $B$, an activation function $\sigma$, and parameters $W$. Let $g(\cdot)$ be a function
  \[
    g(\cdot) = \left(
      \begin{matrix}
        g_1(\cdot)
        \\
        g_2(\cdot)
        \\
        \cdots
        \\
        g_n(\cdot)
        \\
        1
      \end{matrix}
    \right)
  \]
where $g_1(\cdot),\cdots, g_n(\cdot)$ are functions mapping a vector to a number. The operator $f^{B, \sigma}_W$  is defined by
\begin{align}
  \underbrace{f^{B, \sigma}_W \overbrace{(g)}^{\mathclap{\substack{\text{$f_W^{B,\sigma}$'s argument is}\\\text{a function $g$}}}}}_{\mathclap{\substack{\text{$f_W^{B,\sigma}(g)$ is}\\\text{another function}}}} (\cdummy) := \sigma \left( W \left(
  \begin{array}{c}
    \frac{g_1 (\cdot) - \textbf{mean} (g_1, B)}{\sqrt{\textbf{var} (g_1, B)}}\\
    \frac{g_2 (\cdot) - \textbf{mean} (g_2, B)}{\sqrt{\textbf{var} (g_2, B)}}\\
    \vdots\\
    \frac{g_n (\cdot) - \textbf{mean} (g_n, B)}{\sqrt{\textbf{var} (g_n, B)}}\\
    1
  \end{array} \right) \right) \label{eq:zvnklxdfjaljvsaf}
\end{align}
where $\textbf{mean}(r, B)$ is defined by
\[\textbf{mean} (r, B) := \frac{1}{| B |} \sum_{b \in B} r (b).\]
and $\textbf{var}(r, B)$ is defined by
\[ \textbf{var} (r, B) = \textbf{mean} (r^2, B) - \textbf{mean} (r, B)^2 . \]
Note that the first augment of $\textbf{mean}(\cdot, \cdot)$ and
$\textbf{var}(\cdot, \cdot)$ is a function, and the second augment is a set of samples.
A $m$-layer's neural network with BN can be represented by a function
\iftoggle{aistats}{
  \begin{align*}
    F^B_{\{W_i\}_{i=1}^m}(\cdot) :=&
    f_{W_m}^{B, \sigma_m}(f_{W_{m-1}}^{B, \sigma_{m-1}}(\cdots (f_{W_1}^{B, \sigma_1}(I))))(\cdot)\quad \text{or}\\
    & f_{W_m}^{B, \sigma_m}\circ f_{W_{m-1}}^{B, \sigma_{m-1}}\circ\cdots\circ f_{W_1}^{B, \sigma_1}(I)(\cdot)
  \end{align*}
}{
  \begin{align*}
    F^B_{\{W_i\}_{i=1}^m}(\cdot) :=
    f_{W_m}^{B, \sigma_m}(f_{W_{m-1}}^{B, \sigma_{m-1}}(\cdots (f_{W_1}^{B, \sigma_1}(I))))(\cdot)\quad \text{or}\quad f_{W_m}^{B, \sigma_m}\circ f_{W_{m-1}}^{B, \sigma_{m-1}}\circ\cdots\circ f_{W_1}^{B, \sigma_1}(I)(\cdot)
  \end{align*}
}
where $I(\cdot)$ is the identical mapping function from a vector to itself, that
is, $I(x) = x$.

\paragraph{BN's objective is very different from the objective in our common sense}
Using the same way, we can represent a fully connected network in our common sense (without BN) by
\begin{align*}
F_{\{W_i\}_{i=1}^m}(\cdot) := \bar{f}_{W_m}^{\sigma_m}\circ \bar{f}_{W_{m-1}}^{\sigma_{m-1}}\circ\cdots\circ \bar{f}_{W_1}^{\sigma_1}(I)(\cdot)
\end{align*}
where operator $\bar{f}_W^{\sigma}$ is defined by
\begin{align}
\bar{f}^{\sigma}_W(g)(\cdot) := \sigma(Wg(\cdot)).\label{eq:czldjsalvmalksfj}
\end{align}
Besides of the difference of operator function definitions, their ultimate objectives are also different. Given the training date set $\mathcal{D}$, the objective without BN (or the objective in our common sense) is defined by
  \begin{equation}
    \min_{\{W_j\}_{j=1}^m} \quad {1\over |\mathcal{D}|}\sum_{(\x, y)\in \mathcal{D}} l(F_{\{W_j\}_{j=1}^m}(\x), y),\label{eq:zhnklasfuoeq}
  \end{equation}
where $l(\cdot, \cdot)$ is a predefined loss function, while the objective with BN (or the objective BN is actually solving) is
\begin{equation}
({\bf BN})\quad \min_{\{W_j\}_{j=1}^m} \quad {1\over |\mathcal{B}|}\sum_{B\in \mathcal{B}}{1\over |B|} \sum_{(\x, y)\in B} l(F^B_{\{W_j\}_{j=1}^m}(\x), y),\label{eq:vjlkasfnkl}
\end{equation}
where $\mathcal{B}$ is the set of batches.

Therefore, the objective of BN could be very different from the objective in our common sense in general.

\paragraph{BN could be very sensitive to the sample strategy and the minibatch size}
The super set $\mathcal{B}$ has different forms depending on how to define mini-batchs. For example,
  \begin{itemize}
  \item People can choose $b$ as the size of minibatch, and form $\mathcal{B}$ by
  \[\mathcal{B} := \left\{B\subset \mathcal{D}:~|B| = b\right\}. \]
  Choosing $\mathcal{B}$ in this way implicitly assumes that all nodes can
  access the same dataset, which may not be true in practice.
  \item If data is distributed ($\mathcal{D}_i$ is the local dataset on the $i$-th
    worker, disjoint with others, satisfying $\mathcal{D}=\mathcal{D}_1\cup
    \cdots \cup \mathcal{D}_n$), a typical $\mathcal{B}$ is defined by $\mathcal{B} = \mathcal{B}_1 \cup \mathcal{B}_2 \cup \cdots \cup \mathcal{B}_n$ with $\mathcal{B}_i$ defined by
      \[\mathcal{B}_i := \left\{B\subset \mathcal{D}_i:~|B| = b\right\}. \]
  \item When the mini-batch is chosen to be the whole dataset, $\mathcal{B}$ contains only one
    element $\mathcal{D}$.
  \end{itemize}
After figure out the implicit objective BN optimizes in~\eqref{eq:vjlkasfnkl}, it is not difficult to have following key observations
\begin{itemize}
\item The BN objectives vary a lot when different sampling strategies are
  applied. This explains why the convergent solution could be very different when we change the sampling strategy;
  \item For the same sample strategy, BN's objective function also varies if the size of minibatch gets changed. This explains why BN could be sensitive to the batch size.
\end{itemize}
The observations above may seriously bother us how to appropriately choose parameters for BN, since it does not optimize the objective in our common sense, and could be sensitive to the sampling strategy and the size of minibatch.
\paragraph{BN's objective \eqref{eq:vjlkasfnkl} with
  $\mathcal{B}=\{\mathcal{D}\}$ has the same optimal value as the original
  objective \eqref{eq:zhnklasfuoeq}} When the batch size is equal to the size of
the whole dataset in \eqref{eq:vjlkasfnkl}, the objective becomes
\begin{equation}
(\textbf{FN})\quad\min_{\{W_j\}_{j=1}^m} \quad {1\over |\mathcal{D}|} \sum_{(\x, y)\in \mathcal{D}} l(F^{\mathcal{D}}_{\{W_j\}_{j=1}^m}(\x), y),
\label{eq:fn}
\end{equation}
which differs from the original objective \eqref{eq:zhnklasfuoeq} only in the
first argument of $l$. Noting the only difference between \eqref{eq:zvnklxdfjaljvsaf} and
\eqref{eq:czldjsalvmalksfj} when $B$ is a constant $\mathcal{D}$ is a linear
transformation which can be absorbed into $W$:
\setlength{\arraycolsep}{1.5pt}
\begin{align*}
  & f^{\mathcal{D}, \sigma}_W (g) (\cdummy) = \sigma \left( W \left(
  \begin{array}{c}
    \frac{g_1 (\cdummy) - \textbf{mean} (g_1, \mathcal{D})}{\sqrt{\textbf{var}
    (g_1, \mathcal{D})}}\\
    \frac{g_2 (\cdummy) - \textbf{mean} (g_2, \mathcal{D})}{\sqrt{\textbf{var}
    (g_2, \mathcal{D})}}\\
    \vdots\\
    \frac{g_n (\cdummy) - \textbf{mean} (g_n, \mathcal{D})}{\sqrt{\textbf{var}
    (g_n, \mathcal{D})}}\\
    1
  \end{array} \right) \right)\\
  = & \sigma \left(  \underbrace{W
\scriptsize{
     \left( \begin{array}{ccccc}
    \frac{1}{\sqrt{\textbf{var} (g_1, \mathcal{D})}} \hspace{-5mm}&  &  &  & \frac{-
    \textbf{mean} (g_1, \mathcal{D})}{\sqrt{\textbf{var} (g_1, \mathcal{D})}}\\
    & \frac{1}{\sqrt{\textbf{var} (g_2, \mathcal{D})}}\hspace{-5mm} &  &  & \frac{-
    \textbf{mean} (g_2, \mathcal{D})}{\sqrt{\textbf{var} (g_2, \mathcal{D})}}\\
    &  & \ddots\hspace{-5mm} &  & \vdots\\
    &  &  & \frac{1}{\sqrt{\textbf{var} (g_n, \mathcal{D})}} & \frac{-
    \textbf{mean} (g_n, \mathcal{D})}{\sqrt{\textbf{var} (g_n, \mathcal{D})}}\\
    &  &  &  & 1
  \end{array} \right)}}_{=:W'} \left( \begin{array}{c}
    g_1 (\cdummy)\\
    g_2 (\cdummy)\\
    \vdots\\
    g_n (\cdummy)\\
    1
  \end{array} \right) \right)\\
  = & \sigma \left( W' \left( \begin{array}{ccccc}
    g_1 (\cdummy) & g_2 (\cdummy) & \cdots & g_n (\cdot) & 1
  \end{array} \right)^{\top} \right) .
\end{align*} \setlength{\arraycolsep}{5pt}
If we use this $W'$ as our new $W$, \eqref{eq:zhnklasfuoeq} and
\eqref{eq:vjlkasfnkl} has the same form and the two objectives
should have the same optimal
value, which means when $\mathcal{B}=\{\mathcal{D}\}$ adding BN does not hurt
the expressiveness of the network. However, since each layer's input has been
normalized, in general the condition number of the objective is reduced, and
thus easier to train.

\section{Solving full normalization objective \eqref{eq:fn} via compositional optimization}

The BN objective contains the normalization operator which mostly reduces the condition number of the objective that can accelerate the training process. While the FN formulation in \eqref{eq:fn} provides a more stable and trustable way to define the BN formation, it is very challenging to solve the FN formulation in practice. If follow the standard SGD used in solving BN to solve \eqref{eq:fn}, then every single iteration needs to involve the whole dataset since $B = \mathcal{D}$, which is almost impossible in practice. The key difficulty to solve  \eqref{eq:fn} lies on that there exists ``expectation'' in each layer, which is very expensive to compute even we just need to compute a stochastic gradient if use the standard SGD method.

To solve \eqref{eq:fn} efficiently, we follow the spirit of compositional
optimization to develop a new algorithm namely, \textsc{Multilayer Compositional
Stochastic Gradient Descent} (\Cref{algo:theory}). The proposed algorithm does not have any requirement on the size of minibatch. The key idea is to estimate the expectation from the current samples and historical record.

\subsection{Formulation}
To propose our solution to solve~\eqref{eq:fn}, let us define the objective in a more general but neater way in the following.

  When $\mathcal{B}=\{\mathcal{D}\}$, \eqref{eq:vjlkasfnkl} is a special case
  of the following general objective:
  \begin{align}
    \label{eq:4f20e2e09f3b164b9b6851adaa90c5cecba511ad}
    \min_w f(w) := & \mathbb{E}_{\xi} \Big[  F_1^{w_1, \mathcal{D}}  \circ F_2^{w_2, \mathcal{D}} \circ \cdots \circ F_n ^{w_n, \mathcal{D}} \circ I (\xi)\Big] ,
  \end{align}
  where $\xi$ represents a sample in the dataset $\mathcal{D}$, for example it
  can be the pair $(\x,y)$ in \eqref{eq:vjlkasfnkl}. $w_i$ represents the
  parameters of the $i$-th layer. $w$ represents all parameters of the layers:
  $w:=(w_1,\ldots, w_n)$. Each $F_i$ is an operator of the form:
  \begin{align*}
    F_i^{w_i, \mathcal{D}}(g)(\cdot) := & F_i(w_i; g(\cdot); \mathbb{E}_{\xi\in\mathcal{D}}e_i(g(\xi))),
  \end{align*}
  where $e_i$ is used to compute statistics over the dataset. For example, with
  $e_i(x) = [x, x^2]$, the mean and variance of the layer's input over the
  dataset can be calculated, and $F_i^{w_i,\mathcal{D}}$ can use that
  information, for example, to perform normalization.

There exist compositional optimization
algorithms \citep{wang2016stochastic,wang2016accelerating} for
solving \eqref{eq:4f20e2e09f3b164b9b6851adaa90c5cecba511ad} when $n=2$, but for
$n>2$ we still do not have a good algorithm to solve it. We follow the spirit to extend the
compositional optimization algorithms to solve the general optimization problem
\eqref{eq:4f20e2e09f3b164b9b6851adaa90c5cecba511ad} as shown in
\Cref{algo:theory}. See \Cref{sec:exp} for an implementation when it comes to
normalization.

  To simplify notation, given $w=(w_1,\ldots, w_n)$ we define:
  \begin{align*}
   \mathbf{e}_i(w; \xi) :=&  e_i({F}_{i+1}^{w_{i+1}, \mathcal{D}}\circ
    {F}_{i+2}^{w_{i+2}, \mathcal{D}} \circ \cdots \circ
    {F}_{n}^{w_{n}, \mathcal{D}}\circ I (\xi)).
  \end{align*}
  We define the following operator if we already have an estimation, say
  $\tilde{e}_i$, of $\mathbb{E}_{\xi\in\mathcal{D}}\mathbf{e}_i(w; \xi)$:
  \begin{align*}
    \tilde{F}_i^{w_i, \tilde{e}_i}(g)(\cdot) := & F_i(w_i; g(\cdot); \tilde{e}_i).
  \end{align*}
  Given $w=(w_1,\ldots, w_n)$ and $\tilde{e}=(\tilde{e}_1,\ldots, \tilde{e}_n)$, we define
  \iftoggle{aistats}{
    \begin{align*}
      \tilde{\mathbf{e}}_i(w; \xi; \tilde{e}) :=&  e_i(\tilde{F}_{i+1}^{w_{i+1}, \tilde{e}_{i+1}}\circ
                                                  \tilde{F}_{i+2}^{w_{i+2}, \tilde{e}_{i+2}} \circ\\
                                                & \cdots \circ \tilde{F}_{n}^{w_{n}, \tilde{e}_{n}}\circ I (\xi)).
    \end{align*}
  }{
    \begin{align*}
      \tilde{\mathbf{e}}_i(w; \xi; \tilde{e}) :=&  e_i(\tilde{F}_{i+1}^{w_{i+1}, \tilde{e}_{i+1}}\circ
                                                  \tilde{F}_{i+2}^{w_{i+2}, \tilde{e}_{i+2}} \circ \cdots \circ \tilde{F}_{n}^{w_{n}, \tilde{e}_{n}}\circ I (\xi)).
    \end{align*}
  }

\begin{algorithm}
  \caption{Multilayer Compositional
Stochastic Gradient Descent (MCSGD) algorithm\label{algo:theory}}
  \begin{algorithmic}[1]
    \begin{xiangru}
      \Require Learning rate $\{\gamma_k\}_{k=0}^K$, approximation rate
      $\{\alpha_k\}_{k=0}^K$, dataset $\mathcal{D}$, initial point
      $w^{(0)}:=(w_1^{(0)},\ldots, w_n^{(0)})$ and initial estimations
      $\tilde{e}^{(0)}:=(\tilde{e}_1^{(0)},\ldots, \tilde{e}_n^{(0)})$.
      \For {$k=0,1,2,\ldots, K$}
      \State For each $i$, randomly select a sample $\xi_k$, estimate $\mathbb{E}_\xi[\mathbf{e}_i(w^{(k)}; \xi)]$ by \[\tilde{e}_i^{(k+1)} \gets
      (1-\alpha_k)\tilde{e}_i^{(k)} + \alpha_k \tilde{\mathbf{e}}_i(w^{(k)}; \xi_k; \tilde{e}^{(k)}).
      \]
      \State Ask the oracle $\mathcal{O}$ using $w^{(k)}$ and estimated means
      $\tilde{e}^{(k+1)}$ to obtain the approximated gradient at $w^{(k)}$:
      $g^{(k)}$.

      \Comment{See \Cref{remark:oracle} for discussion on the oracle.}

      \State $w^{(k+1)}\gets w^{(k)} - \gamma_k g^{(k)}$.
      \EndFor
    \end{xiangru}
  \end{algorithmic}
\end{algorithm}
\begin{remark}[MCSGD oracle]
  \label{remark:oracle}
  In MCSGD, the gradient oracle takes the current parameters of the
  model, and the current estimation of each
  $\mathbb{E}_\xi[\mathbf{e}_i(w^{(k)}; \xi)]$, to output an estimation of a
  stochastic gradient.

  For example for an objective like
  \eqref{eq:4f20e2e09f3b164b9b6851adaa90c5cecba511ad}, the derivative of the
  loss function w.r.t. the $i$-th layer's parameters is
  \iftoggle{aistats}{
    \begin{align*}
      &\partial_{w_i} f (w) \\
      = & \partial_{w_i} (\mathbb{E}_{\xi} [F_1^{w_1,
          \mathcal{D}} \circ F_2^{w_2, \mathcal{D}} \circ \cdots \circ F_n^{w_n,
          \mathcal{D}} \circ I (\xi)])\\
      = & \mathbb{E}_{\xi} \left[ \begin{array}{l}
                                    {}[\partial_x (F_1^{w_1, \mathcal{D}} (x) (\xi))]_{x = F_2^{w_2,
                                    \mathcal{D}} \circ \cdots \circ F_n^{w_n, \mathcal{D}} \circ I} \cdot\\
                                    {}[\partial_x (F_2^{w_2, \mathcal{D}} (x) (\xi))]_{x = F_3^{w_3,
                                    \mathcal{D}} \circ \cdots \circ F_n^{w_n, \mathcal{D}} \circ I} \cdot
                                    \cdots \cdot\\
                                    {}[\partial_x (F_{i - 1}^{w_{i - 1}, \mathcal{D}} (x) (\xi))]_{x =
                                    F_i^{w_i, \mathcal{D}} \circ \cdots \circ F_n^{w_n, \mathcal{D}} \circ I}
                                    \cdot\\
                                    {}[\partial_{w_i} (F_i^{w_i, \mathcal{D}} (x) (\xi))]_{x = F_i^{w_{i + 1},
                                    \mathcal{D}} \circ \cdots \circ F_n^{w_n, \mathcal{D}} \circ I}
                                  \end{array} \right], \numberthis\label{eq:vnklastnwkqr}
    \end{align*}
  }{
    \begin{align*}
      \partial_{w_i} f (w)
      = & \partial_{w_i} (\mathbb{E}_{\xi} [F_1^{w_1,
          \mathcal{D}} \circ F_2^{w_2, \mathcal{D}} \circ \cdots \circ F_n^{w_n,
          \mathcal{D}} \circ I (\xi)])\\
      = & \mathbb{E}_{\xi} \left[ \begin{array}{l}
                                    {}[\partial_x (F_1^{w_1, \mathcal{D}} (x) (\xi))]_{x = F_2^{w_2,
                                    \mathcal{D}} \circ \cdots \circ F_n^{w_n, \mathcal{D}} \circ I} \cdot\\
                                    {}[\partial_x (F_2^{w_2, \mathcal{D}} (x) (\xi))]_{x = F_3^{w_3,
                                    \mathcal{D}} \circ \cdots \circ F_n^{w_n, \mathcal{D}} \circ I} \cdot
                                    \cdots \cdot\\
                                    {}[\partial_x (F_{i - 1}^{w_{i - 1}, \mathcal{D}} (x) (\xi))]_{x =
                                    F_i^{w_i, \mathcal{D}} \circ \cdots \circ F_n^{w_n, \mathcal{D}} \circ I}
                                    \cdot\\
                                    {}[\partial_{w_i} (F_i^{w_i, \mathcal{D}} (x) (\xi))]_{x = F_i^{w_{i + 1},
                                    \mathcal{D}} \circ \cdots \circ F_n^{w_n, \mathcal{D}} \circ I}
                                  \end{array} \right], \numberthis\label{eq:vnklastnwkqr}
    \end{align*}
  }
  where for any $j \in [1, \ldots, i - 1]$:
  \iftoggle{aistats}{
\begin{align*}
  & [\partial_x F_j^{w_j, \mathcal{D}} (x) (\xi)]_{x = F_{j + 1}^{w_{j + 1},
  \mathcal{D}} \circ \cdots \circ F_n^{w_n, \mathcal{D}} \circ I}\\
  = & \left[ \begin{array}{c}
    \partial_x F_j (w_j; x ; y) +\\
    \partial_y F_j (w_j; x ; y) \cdot \mathbb{E}_{\xi'\sim \mathcal{D}} D_j (\xi')
             \end{array} \right]\\
  & \text{with}{\scriptsize{\begin{array}{l}
    x = F_{j + 1}^{w_{j + 1}, \mathcal{D}} \circ \cdots \circ F_n^{w_n,
    \mathcal{D}} \circ I (\xi),\\
    y =\mathbb{E}_{\xi} [\mathbf{e}_j (w ; \xi)],\\
    D_j (\xi') = [\partial_z e_i (z)]_{z = F_{j + 1}^{w_{j + 1}, \mathcal{D}}
    \circ \cdots \circ F_n^{w_n, \mathcal{D}} \circ I (\xi')}
  \end{array}}}.
\end{align*}
}{
\begin{align*}
  {}[\partial_x F_j^{w_j, \mathcal{D}} (x) (\xi)]_{x = F_{j + 1}^{w_{j + 1},
  \mathcal{D}} \circ \cdots \circ F_n^{w_n, \mathcal{D}} \circ I} = & \left[
  \begin{array}{c}
    \partial_x F_j (w_j ; x ; y) +\\
    \partial_y F_j (w_j ; x ; y) \cdot \mathbb{E}_{\xi' \sim \mathcal{D}} D_j
    (\xi')
  \end{array} \right]\\
  & \text{with} \left\{\begin{array}{l}
    x = F_{j + 1}^{w_{j + 1}, \mathcal{D}} \circ \cdots \circ F_n^{w_n,
    \mathcal{D}} \circ I (\xi)\\
    y =\mathbb{E}_{\xi} [\mathbf{e}_j (w ; \xi)]\\
    D_j (\xi') = [\partial_z e_i (z)]_{z = F_{j + 1}^{w_{j + 1}, \mathcal{D}}
    \circ \cdots \circ F_n^{w_n, \mathcal{D}} \circ I (\xi')}
  \end{array}\right. .
\end{align*}
}
The oracle samples a $\xi'$ for each layer and calculate this derivative with
$y$ set to an estimated value $\hat{e}_j$, and returns the stochastic gradient.
The expectation of this stochastic gradient will be {\eqref{eq:vnklastnwkqr}} with some error caused by the difference between the
estimated $\hat{e}_j$ and the true expectation $\mathbb{E}_{\xi}
[\mathbf{e}_j (w ; \xi)]$. The more accurate the estimation is, the closer
the expectation of this stochastic and the true gradient {\eqref{eq:vnklastnwkqr}} will be.

In practice, we often use the same $\xi'$ for each layer to save some
computation resource, this introduce some bias the expectation of the
estimated stochastic gradient. One such implementation can be found in \Cref{algo:impl}.
\end{remark}

\subsection{Theoretical analysis}
In this section we analyze \Cref{algo:theory} to show its convergence when
applied on \eqref{eq:4f20e2e09f3b164b9b6851adaa90c5cecba511ad} (a generic version of the FN formulation in \eqref{eq:fn}). The detailed
proof is provided in the supplementary material. We use the following assumptions as
shown in \Cref{ass:a7f57b116c8ed2da7490c29a079fe3e51df929cc} for the analysis.
\begin{assumption}\label{ass:a7f57b116c8ed2da7490c29a079fe3e51df929cc}
  \begin{enumerate}
  \item The gradients $g^{(k)}$'s are bounded:
    \begin{equation*}
      \|g^{(k)}\| \leqslant \mathscr{G}, \forall k
    \end{equation*}
    for some constant $\mathscr{G}$.

  \item The variance of all $e_i$'s are bounded:
      \begin{equation*}
        \mathbb{E}_{\xi}\|\mathbf{e}_i(w;\xi) - \mathbb{E}_{\xi}\mathbf{e}_i(w;\xi)\|^2 \leqslant \sigma^2,\forall i, w,
      \end{equation*}
      \begin{equation*}
        \mathbb{E}_{\xi}\|\mathbf{e}_i(w;\xi; \hat{e}) - \mathbb{E}_{\xi}\mathbf{e}_i(w;\xi;\hat{e})\|^2 \leqslant \sigma^2,\forall i , w,\hat{e}
      \end{equation*}
    for some constant $\sigma$.

  \item The error of approximated gradient $\mathbb{E}[g^{(k)}]$ is proportional to the
    error of approximation for $\mathbb{E} [e_i]$:
      \iftoggle{aistats}{
        \begin{align*}
          & \|\mathbb{E}_{\xi_k}[g^{(k)}] - \nabla f(w^{(k)})\|^2 \\ \leqslant & L_g \sum_{i=1}^{n}\| \tilde{e}_i^{(k+1)} - \mathbb{E}_{\xi_k} [\mathbf{e}_i(w^{(k)};\xi_k)] \|^2, \forall k
        \end{align*}
      }{
        \begin{equation*}
          \|\mathbb{E}_{\xi_k}[g^{(k)}] - \nabla f(w^{(k)})\|^2 \leqslant L_g \sum_{i=1}^{n}\| \tilde{e}_i^{(k+1)} - \mathbb{E}_{\xi_k} [\mathbf{e}_i(w^{(k)};\xi_k)] \|^2, \forall k
        \end{equation*}
      }
    for some constant $L_g$.

  \item  All functions and their first order derivatives are Lipschitzian with Lipschitz constant $L$.

  \item  The minimum of the objective $f(w)$ is finite.

  \item $\gamma_k,\alpha_k$ are monotonically decreasing. $\gamma_k =
    O(k^{-\gamma})$ and $\alpha_k = O(k^{-a})$ for some constants $\gamma > a > 0$.

  \end{enumerate}
\end{assumption}

It can be shown under the given assumptions, the approximation errors will vanish shown in \Cref{lemma:98620d3402422d993259e3a830cd49a398b44a20}.

\begin{lemma}[Approximation error]
  \label{lemma:98620d3402422d993259e3a830cd49a398b44a20}
  Choose the learning rate $\gamma_k$ and $\alpha_k$ in Algorithm~\ref{algo:theory} in the form defined in Assumption~\ref{ass:a7f57b116c8ed2da7490c29a079fe3e51df929cc}-6 with parameters $\gamma$ and $a$. Under \Cref{ass:a7f57b116c8ed2da7490c29a079fe3e51df929cc}, the sequence generated in Algorithm~\ref{algo:theory} satisfies
    \iftoggle{aistats}{
      \begin{align*}
        &\mathbb{E} \| \tilde{e}_{i}^{(k + 1)} -\mathbb{E}_{\xi} \mathbf{e}_{i} (w^{(k)} ; \xi) \|^2 \\ \leqslant &\mathscr{E}(k^{- 2 \gamma + 2 a + \varepsilon} + k^{- a + \varepsilon}), \forall i, \forall k, \numberthis \label{eq:b5cf96af2f633135e356bd345b0fa5731e5e2b15}
      \end{align*}
      \begin{align*}
        & \mathbb{E} \| \tilde{e}_i^{(k + 1)} -\mathbb{E}_{\xi} \mathbf{e}_i (w^{(k)} ; \xi)
        \|^2\\ \leqslant & \left( 1 - \frac{\alpha_k}{2} \right) \mathbb{E} \|
                         \tilde{e}_i^{(k)} -\mathbb{E}_{\xi} \mathbf{e}_i (w^{(k - 1)} ; \xi) \|^2\\
                       & +\mathscr{C} (k^{- 2 \gamma + a + \varepsilon} + k^{- 2 a +
                         \varepsilon}), \quad \forall i, \forall k.\numberthis \label{eq:2ab3f273dad20de3dbd553953676c02229427496}
      \end{align*}
    }{
      \begin{align*}
        \mathbb{E} \| \tilde{e}_{i}^{(k + 1)} -\mathbb{E}_{\xi} \mathbf{e}_{i} (w^{(k)} ; \xi) \|^2 \leqslant & \mathscr{E}(k^{- 2 \gamma + 2 a + \varepsilon} + k^{- a + \varepsilon}), \forall i, \forall k. \numberthis \label{eq:b5cf96af2f633135e356bd345b0fa5731e5e2b15}
\\
        \mathbb{E} \| \tilde{e}_i^{(k + 1)} -\mathbb{E}_{\xi} \mathbf{e}_i (w^{(k)} ; \xi)
        \|^2 \leqslant & \left( 1 - \frac{\alpha_k}{2} \right) \mathbb{E} \|
                         \tilde{e}_i^{(k)} -\mathbb{E}_{\xi} \mathbf{e}_i (w^{(k - 1)} ; \xi) \|^2\\
                       & +\mathscr{C} (k^{- 2 \gamma + a + \varepsilon} + k^{- 2 a +
                         \varepsilon}), \quad \forall i, \forall k.\numberthis \label{eq:2ab3f273dad20de3dbd553953676c02229427496}
      \end{align*}
    }
    for any $\varepsilon$ satisfying $1-a > \varepsilon > 0$, where
    $\mathscr{E}$ and $\mathscr{C}$ are two constants independent of $k$.
\end{lemma}

Then it can be shown that on \eqref{eq:4f20e2e09f3b164b9b6851adaa90c5cecba511ad}
\Cref{algo:theory} converges as seen in
\Cref{thm:384575d2d125c25265b5d7bf9c978acb95307110} and
\Cref{coro:8dce0d0cfc00f8fa3138e18c07037ccfc2ffd86a}. It is worth noting that the
convergence rate in \Cref{coro:8dce0d0cfc00f8fa3138e18c07037ccfc2ffd86a} is
slower than the $\frac{1}{\sqrt{K}}$ convergence rate of SGD without any
normalization. This is due to the estimation error. If the estimation
error is small (for example, the samples in a batch are randomly selected and
the batch size is large), the convergence will be fast. %

\begin{theorem}[Convergence]
  \label{thm:384575d2d125c25265b5d7bf9c978acb95307110}
  Choose the learning rate $\gamma_k$ and $\alpha_k$ in
  Algorithm~\ref{algo:theory} in the form defined in
  Assumption~\ref{ass:a7f57b116c8ed2da7490c29a079fe3e51df929cc}-6 with
  parameters $\gamma$ and $a$ satisfying $a< 2\gamma -1, a < 1/2$, and
  $\frac{\gamma_k L_g}{\alpha_{k+1}} \leqslant \frac{1}{2}$. Under
  \Cref{ass:a7f57b116c8ed2da7490c29a079fe3e51df929cc}, for any integer $K>1$ the
  sequence generated by Algorithm~\ref{algo:theory}
  satisfies%
    \begin{align*}
      \frac{\sum_{k = 0}^K \gamma_k \mathbb{E} \| \partial f (w^{(k)})
      \|^2}{\sum_{k = 0}^K \gamma_k} \leqslant & \frac{\mathscr{H}}{\sum_{k = 0}^K \gamma_k},
    \end{align*}
  where $\mathscr{H}$ is a constant independent of $K$.
\end{theorem}
We next specify the choice of $\gamma$ and $a$ to give a more clear convergence rate for our MCSGD algorithm in \Cref{coro:8dce0d0cfc00f8fa3138e18c07037ccfc2ffd86a}.
\begin{corollary}
  \label{coro:8dce0d0cfc00f8fa3138e18c07037ccfc2ffd86a}
  Choose the learning rate $\gamma_k$ and $\alpha_k$ in
  Algorithm~\ref{algo:theory} in the form defined in
  Assumption~\ref{ass:a7f57b116c8ed2da7490c29a079fe3e51df929cc}-6, more specifically, $\gamma_k =
  \frac{1}{2 L_g}(k + 2)^{- 4/5}$ and $\alpha_k = (k + 1)^{- 2/5}$. Under
  \Cref{ass:a7f57b116c8ed2da7490c29a079fe3e51df929cc}, for any integer $K>1$ the
  sequence generated in Algorithm~\ref{algo:theory} satisfies
  \begin{equation*}
    \frac{\sum_{k = 0}^K \mathbb{E} \| \partial f (w^{(k)}) \|^2}{K + 2} \leqslant \frac{\mathscr{H}}{(K + 2)^{1 / 5}},
  \end{equation*}
  where $\mathscr{H}$ is a constant independent of $K$.
\end{corollary}
\begin{figure}
  \centering
  \includegraphics[width=0.48\columnwidth]{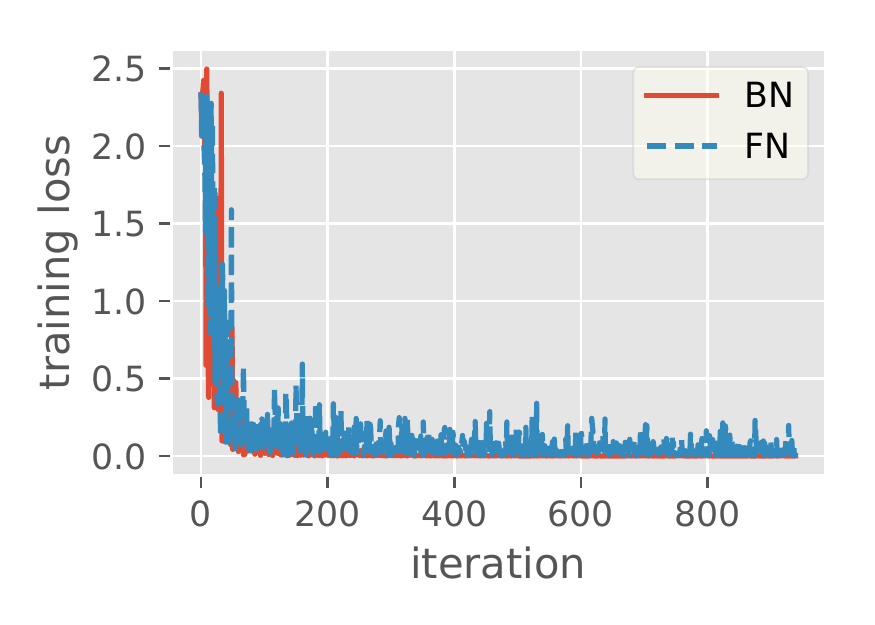}
  \includegraphics[width=0.48\columnwidth]{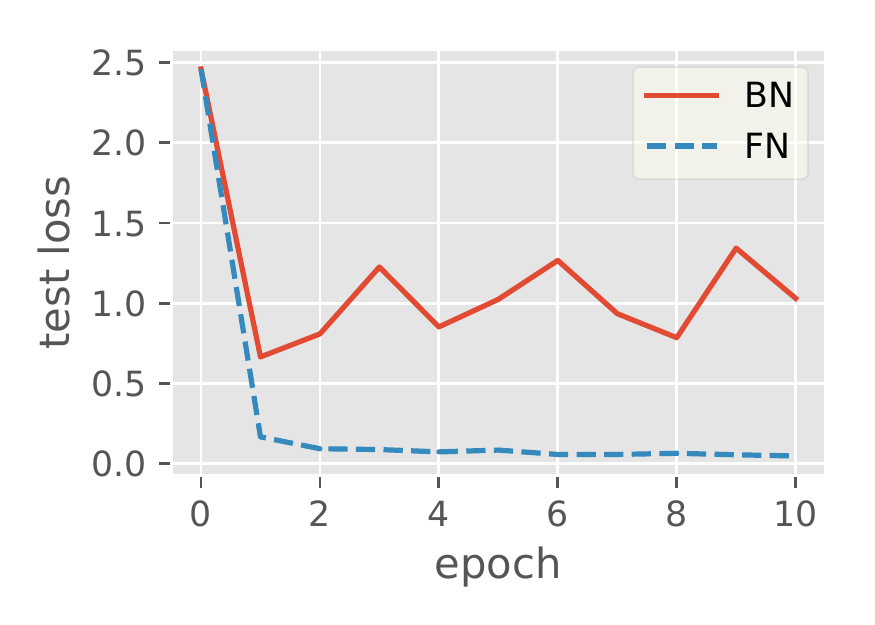}
  \caption{The training and testing error for the given model trained on MNIST
    with batch size 64, learning rate 0.01 and momentum 0.5 using BN or FN for
    the case where all samples in a batch are of a single label. The
    approximation rate $\alpha$ in FN is $(\frac{k}{20} + 1)^{-0.4}$ where $k$
    is the iteration number.}
  \label{fig:conv-mnist-unshuffle}
\end{figure}

\begin{figure}
  \centering
  \includegraphics[width=0.48\columnwidth]{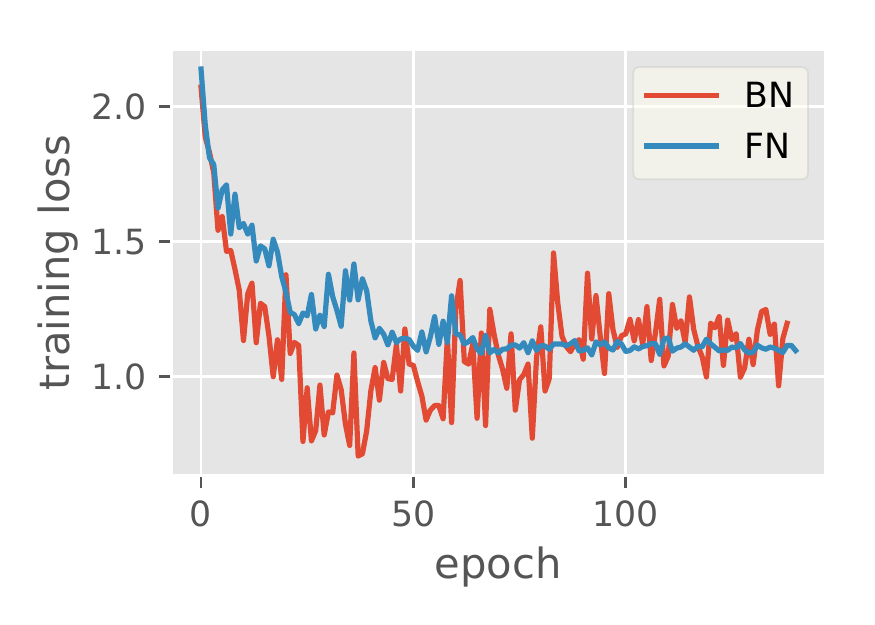}
  \includegraphics[width=0.48\columnwidth]{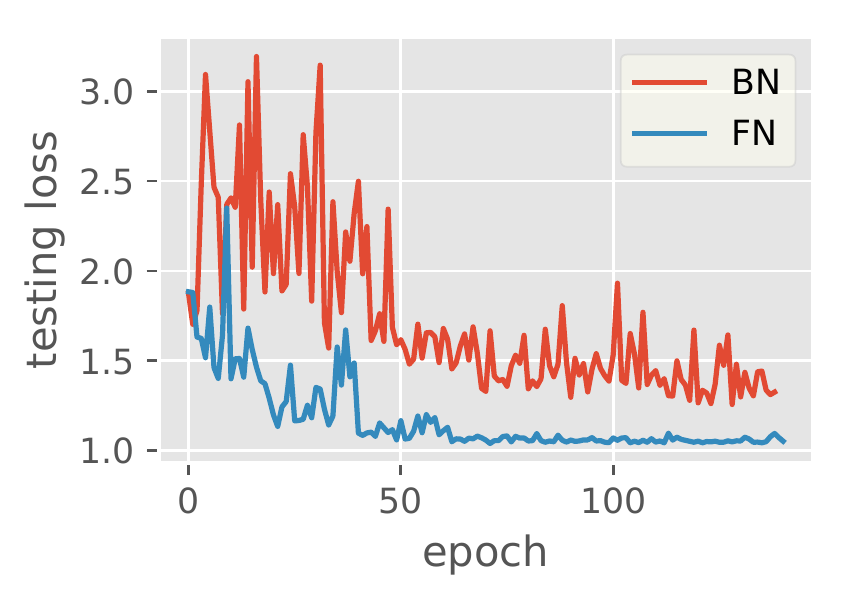}
  \caption{The training and testing error for the given model trained on CIFAR10
    with batch size 64, learning rate 0.01 and momentum 0.9 using BN and FN for the case where all samples in a batch are of no more than 3 labels.
    The learning rate is decreased by a factor of 5 for every 20 epochs. The
    approximation rate $\alpha$ in FN is $(\frac{k}{5} + 1)^{-0.3}$
    where $k$ is the iteration number.}
  \label{fig:conv-cifar-unshuffle}
\end{figure}

\begin{figure*}
  \centering
  \includegraphics[width=\textwidth]{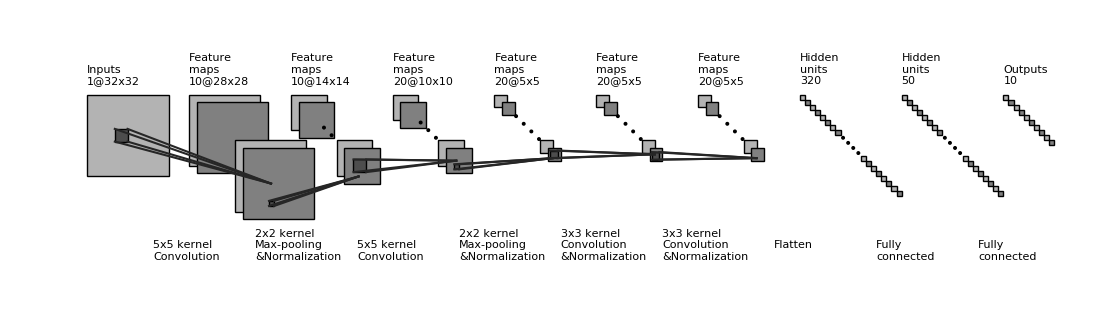}
  \caption{The simple neural network used in the experiments. The
    ``Normalization'' between layers can be BN or FN. The activation
    function is relu and the loss function is negative log likelihood.}
  \label{fig:simplenn}
\end{figure*}

\begin{figure}
  \centering
  \includegraphics[width=0.48\columnwidth]{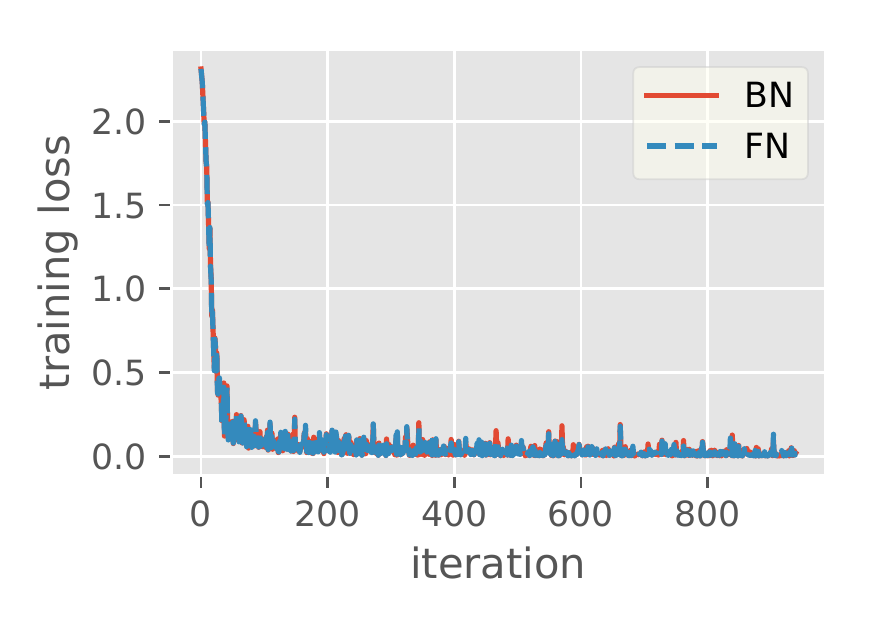}
  \includegraphics[width=0.48\columnwidth]{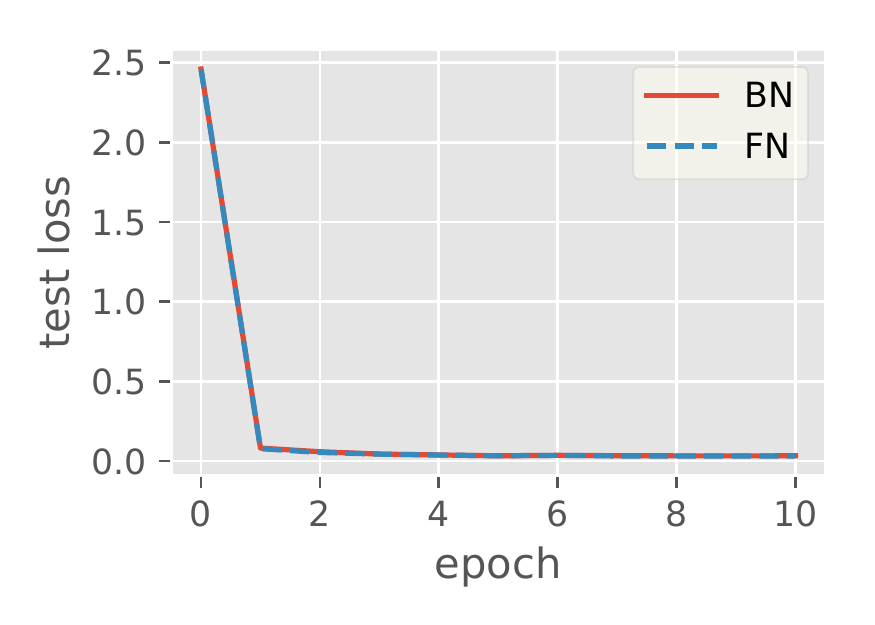}
  \caption{The training and testing error for the given model trained on MNIST
    with batch size 64, learning rate 0.01 and momentum 0.5 using BN or FN for
    the case where samples in a batch are randomly selected. The approximation
    rate $\alpha$ in FN is $(\frac{k}{20} + 1)^{-0.4}$ where $k$ is the
    iteration number.}
  \label{fig:conv-mnist-shuffle}
\end{figure}

\begin{figure}
  \centering
  \includegraphics[width=0.48\columnwidth]{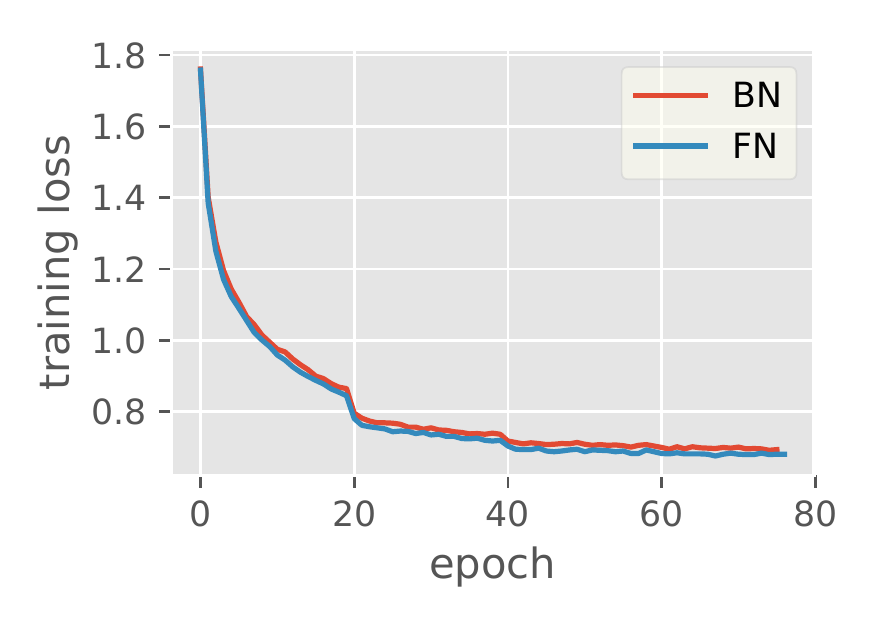}
  \includegraphics[width=0.48\columnwidth]{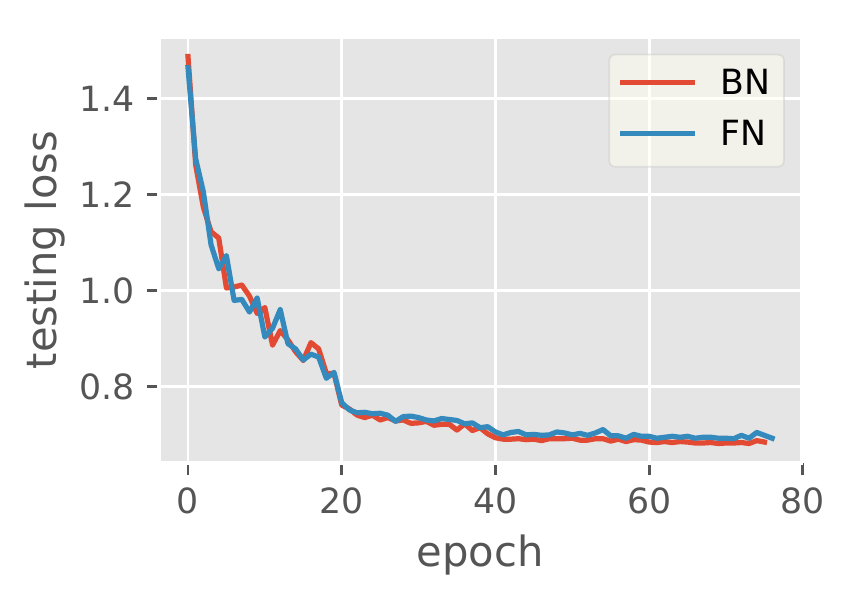}
  \caption{The training and testing error for the given model trained on CIFAR10
    with batch size 64, learning rate 0.01 and momentum 0.9 using BN or FN for
    the case where samples in a batch are randomly selected. The learning rate
    is decreased by a factor of 5 for every 20 epochs. The approximation rate
    $\alpha$ in FN is $(\frac{k}{20} + 1)^{-0.2}$ where $k$ is the iteration
    number.}
  \label{fig:conv-cifar-shuffle}
\end{figure}

\section{Experiments}
\label{sec:exp}

Experiments are conducted to validate the effectiveness of our MCSGD algorithm for solving the FN formulation. We consider two settings. The first one in
\Cref{sec:zhlajsfklasjf} shows BN's convergence highly depends on the
size of batches, while FN is more robust to different batch sizes.
The second one in \Cref{sec:zhlajsfklf} shows that FN is more robust to different construction of mini-batches.

We use a simple neural network as the testing network whose architecture is shown in \Cref{fig:simplenn}. Steps 2 and 3 in \Cref{algo:theory} are implemented in \Cref{algo:impl}. The \emph{forward pass} essentially performs Step 2 of \Cref{algo:theory}, which
estimates the mean and variance of layer inputs over the dataset. The
\emph{backward pass} essentially performs Step 3 of \Cref{algo:theory},
which gives the approximated gradient based on current network parameters and
the estimations. Note that as discussed in \Cref{remark:oracle}, for efficiency,
in each iteration of this implementation we are using the same samples to do
estimation in all normalization layers, which saves computational cost but
brings additional bias.

\begin{algorithm}
  \caption{Full Normalization (FN) layer \label{algo:impl}}
  \centering{\fbox{\textsc{Forward pass}}}
  \begin{algorithmic}[1]
    \Require Learning rate $\gamma$, approximation rate $\alpha$, input
    $\mathbf{B}_{\text{in}}\in \mathbb{R}^{b\times d}$, mean estimation $\mu$, and mean
    of square estimation $\nu$. \Comment{$b$ is the batch size, and $d$ is the
      dimension of each sample.}

    \If {training}
      \begin{align*}
        \mu \gets & (1-\alpha) \mu + \alpha \textbf{\textit{mean}}(\mathbf{B}_{\text{in}}), \\
        \nu \gets & (1-\alpha) \nu + \alpha \textbf{\textit{mean\_of\_square}}(\mathbf{B}_{\text{in}}).
      \end{align*}%
    \EndIf

    \State \Return layer
    output: \[\mathbf{B}_{\text{output}} \gets \frac{\mathbf{B}_{\text{in}} - \mu}{\sqrt{\max\{\nu
          - \mu^2, \epsilon\}}}.\] \Comment{$\epsilon$ is a small constant for
      numerical stability.}
  \end{algorithmic}
  \centering{\fbox{\textsc{Backward pass}}}
  \begin{algorithmic}[1]
    \Require Mean estimation $\mu$, mean of square estimation $\nu$, the
    gradient at the output $g_{\text{out}}$, and input $\mathbf{B}_{\text{in}}\in
    \mathbb{R}^{b\times d}$.

    \State Define \[f(\mu, \nu, \mathbf{B}_{\text{in}}) := \frac{\mathbf{B}_{\text{in}} -
      \mu}{\sqrt{\max\{\nu - \mu^2, \epsilon\}}}.\]

    \State \Return the gradient at the input: \[g_{\text{in}} \gets
    g_{\text{out}} \cdot \left(\partial_{\mathbf{B}_{\text{in}}} f + \frac{\partial_{\mu} f + 2
      \mathbf{B}_{\text{in}} \partial_{\nu} f}{bd} \right).\] \Comment{$\partial_a f$ is the
      derivative of $f$ w.r.t. $a$ at $\mu, \nu, \mathbf{B}_{\text{in}}$.}
  \end{algorithmic}
\end{algorithm}

\subsection{Dependence on batch size}
\label{sec:zhlajsfklasjf}

\begin{figure}
  \centering
  \begin{minipage}{0.48\columnwidth}
    \includegraphics[width=\columnwidth]{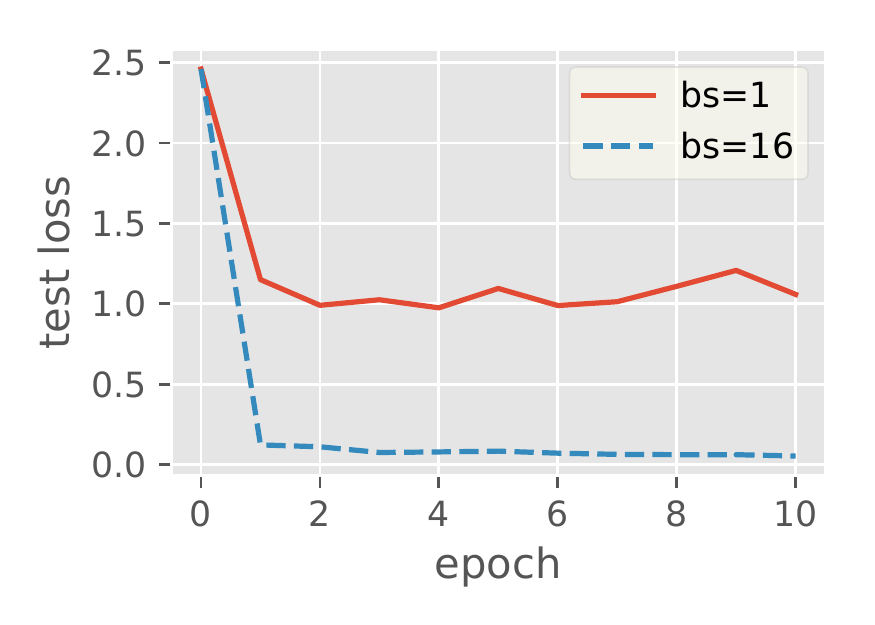}
    \subcaption{With BN layers.}
  \end{minipage}
  \begin{minipage}{0.48\columnwidth}
    \includegraphics[width=\columnwidth]{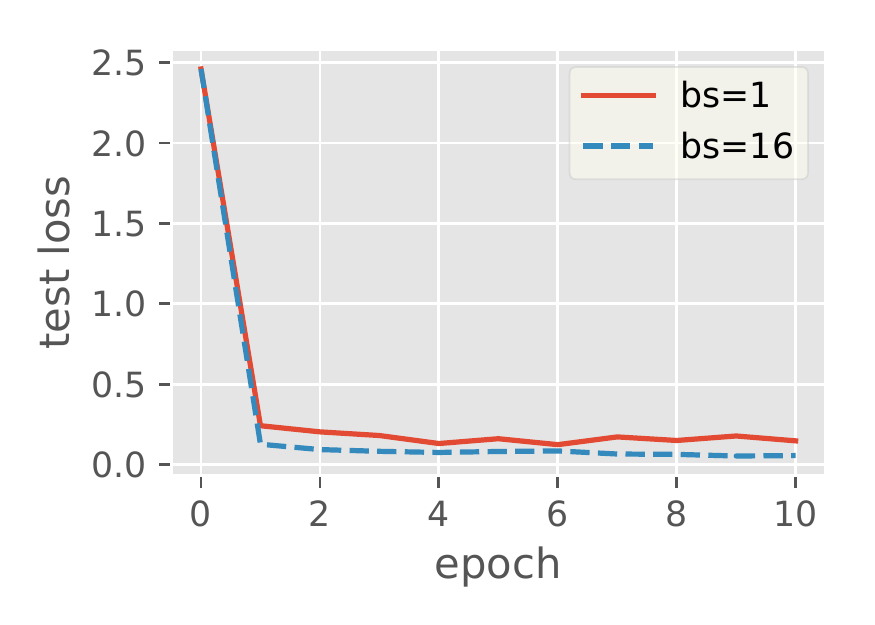}
    \subcaption{With FN layers.}
  \end{minipage}

  \caption{The testing error for the given model trained on MNIST, where each
    sample is multiplied by a random number in $(-2.5, 2.5)$. The optimization
    algorithm is SGD with learning rate is 0.01 and momentum is 0.5. FN gives
    more consistent results under different batch sizes. The approximation rate
    $\alpha$ in FN is $(\frac{k}{20} + 1)^{-0.4}$ where $k$ is the iteration
    number.}
  \label{fig:conv-mnist-batchsize}
\end{figure}

MNIST is used as the testing dataset with some modification by multiplying a random number in $(-2.5, 2.5)$ to each sample to make them more diverse. In this case,
as shown in \Cref{fig:conv-mnist-batchsize}, with batch size 1 and batch size
16, we can see the convergent points are different in BN, while the convergent results of FN are much more close for different batch sizes. Therefore, FN is more robust to the size of mini-batch.

\subsection{Dependence on batch construction}
\label{sec:zhlajsfklf}

We study two cases --- shuffled case and
unshuffled case. For the shuffled case, where the samples in a batch are
randomly sampled, we expect the performance of FN matches BN's in this case. For
the unshuffled case, where the batch contains only samples in just a few number of categories (so the mean and variance are very different from batch to batch). In this case we can observe that the FN outperforms BN.

\paragraph{Unshuffled case}

We show the FN has advantages over BN when the data variation is large among
mini-batches. In the unshuffled setting, we do not shuffle the dataset. It means
that the samples in the same mini-batch are mostly with the same labels. The
comparison uses two datasets (MINST and CIFAR10):
\begin{itemize}
\item On MNIST, the batch size is chosen to be 64 and each batch only contains a single
label. The convergence results are shown in \Cref{fig:conv-mnist-unshuffle}. We
can observe from these figures that for the training loss, BN and FN converge
equally fast (BN might be slightly faster). However since the statistics on any
batch are very different from the whole dataset, the estimated mean and variance
in BN are very different from the true mean and variance on the whole dataset,
resulting in a very high test loss for BN.
\item On CIFAR10, we observe similar results as shown in
\Cref{fig:conv-cifar-unshuffle}. In this case, we restrict the number of labels
in every batch to be no more than 3. Thus BN's performance on CIFAR10 is slightly better than on MNIST. We see the convergence efficiency of both methods is still
comparable in term of the training loss. However, the testing error for BN is
still far behind the FN.
\end{itemize}

\paragraph{Shuffled case}
For the shuffled case, where the samples in each batch are selected randomly, we
expect BN to have similar performance as FN in this case, since the statistics
in a batch is close to the statistics on the whole dataset. The results for
MNIST are shown in \Cref{fig:conv-mnist-shuffle} and the results for CIFAR10 are
shown in \Cref{fig:conv-cifar-shuffle}. We observe the convergence curves of BN
and FN for both training and testing loss match well.
\section{Conclusion}

We provide new understanding for BN from an optimization perspective by
rigorously defining the optimization objective for BN. BN essentially optimizes
an objective different from the one in our common sense. The implicitly targeted
objective by BN depends on the sampling strategy as well as the minibatch size.
That explains why BN becomes unstable and sensitive in some scenarios. The
stablest objective of BN (called FN formulation) is to use the full dataset as
the mini-batch, since it has equivalence to the objective in our common sense.
But it is very challenging to solve such formulation. To solve the FN objective,
we follow the spirit of compositional optimization to develop MCSGD algorithm to
solve it efficiently. Experiments are also conduct to validate the proposed
method.

\bibliographystyle{abbrvnat}

\newpage
\appendix
\begin{center} {\bf\Large Supplementary Materials}
\end{center}
\section{Proofs}

\begin{lemma}\label{lemma:d083ec1ffe10dc6de84ee1d15cf12c068d8b28ab}
If
\begin{align}
  r_{k + 1} \leqslant & (1 - C_1 k^{- a}) r_k + C_2 (k^{- 2 a + \varepsilon} +
  k^{- 2 \gamma + a + \varepsilon}),\label{eq:popojgnlkzxdf}
\end{align}

for any constant $1 > C_1 > 0, C_2 > 0, 1>a>0, 1-a> \varepsilon > 0$, there exists a constant $\mathscr{E}$ such
that
\begin{align*}
  r_{k + 1} \leqslant & (k^{- 2 \gamma + 2 a + \varepsilon} + k^{- a +
  \varepsilon}) \mathscr{E}.
\end{align*}
\end{lemma}

\paragraph{Proof to \Cref{lemma:d083ec1ffe10dc6de84ee1d15cf12c068d8b28ab}}

First we expand the recursion relation \eqref{eq:popojgnlkzxdf} till $k=1$:

\begin{align*}
  r_{k + 1} \leqslant & (1 - C_1 k^{- a}) r_k + C_2  (k^{- 2 a + \varepsilon}
  + k^{- 2 \gamma + a + \varepsilon})\\
  \leqslant & \left( \prod_{\varkappa = 1}^k (1 - C_1 \varkappa^{- a}) \right)
  r_1 + C_2 \left( \sum_{\varkappa = 1}^k (\varkappa^{- 2 a + \varepsilon} +
  \varkappa^{- 2 \gamma + a + \varepsilon}) \left(  \prod_{\mathcal{k}= \kappa
  + 1}^k (1 - C_1 \mathcal{k}^{- a}) \right) \right). \numberthis \label{eq:czljasdkfja}
\end{align*}

As the first step to bound \eqref{eq:czljasdkfja}, we use the monotonicity of the $\ln(\cdot)$
function to bound the form $\prod_{\mathcal{k}= m}^k (1 - C_1 \mathcal{k}^{-
  a})$ for some $m$. Notice that
\begin{align*}
  \ln \left( \prod_{\mathcal{k}= m}^k (1 - C_1 \mathcal{k}^{- a}) \right) = &
  \sum_{\mathcal{k}= m}^k \ln (1 - C_1 \mathcal{k}^{- a})
  \leqslant - C_1 \sum_{\mathcal{k}= m}^k \mathcal{k}^{- a}.
\end{align*}
Then it follows from the monotonicity of $\ln(\cdot)$ that:
\begin{align*}
\prod_{\mathcal{k}= m}^k (1 - C_1 \mathcal{k}^{- a})
  \leqslant & \exp \left( - C_1 \sum_{\mathcal{k}= m}^k \mathcal{k}^{- a}
  \right)
  \leqslant \exp (- C_1  (k - m + 1) k^{- a}) . \numberthis \label{eq:znlkdjasnvla}
\end{align*}

Using this inequality \eqref{eq:czljasdkfja} can be bounded by the following inequality:
\begin{align*}
  r_{k + 1} \leqslant & \left( \prod_{\varkappa = 1}^k (1 - C_1 \varkappa^{-
  a}) \right) r_1 + C_2 \left( \sum_{\varkappa = 1}^k (\varkappa^{- 2 a + \varepsilon} +
  \varkappa^{- 2 \gamma + a + \varepsilon}) \left(  \prod_{\mathcal{k}= \kappa
  + 1}^k (1 - C_1 \mathcal{k}^{- a}) \right) \right)\\
  \overset{\eqref{eq:znlkdjasnvla}}{\leqslant} & \exp (- k^{1 - a}) r_1 + \underbrace{C_2 \left( \sum_{\varkappa = 1}^k (\varkappa^{- 2 a +
  \varepsilon} + \varkappa^{- 2 \gamma + a + \varepsilon}) \exp (- C_1  (k -
  \varkappa) k^{- a}) \right)}_{=:T_k} .\numberthis \label{eq:vnkslkfjblkaf}
\end{align*}
The final step will be bounding $T_k$ above. The term can actually be bounded by a power
of reciprocal function. Notice that for any $1 - a > \varepsilon > 0$, we have the following inequality:
\begin{align*}
  T_k = & C_2  \sum_{\varkappa = 1}^k (\varkappa^{- 2 a +
  \varepsilon} + \varkappa^{- 2 \gamma + a + \varepsilon}) \exp (- C_1  (k -
  \varkappa) k^{- a})\\
  \leqslant & C_2  \sum_{\varkappa = 1}^{k - \lfloor k^{a + \varepsilon}
  \rfloor} (\varkappa^{- 2 a + \varepsilon} + \varkappa^{- 2 \gamma + a +
  \varepsilon}) \exp (- C_1 \lfloor k^{a + \varepsilon} \rfloor k^{- a})
   + C_2  \sum_{\varkappa = k - \lfloor k^{a + \varepsilon} \rfloor + 1}^k
  (\varkappa^{- 2 a + \varepsilon} + \varkappa^{- 2 \gamma + a +
  \varepsilon})\\
  = & O (\exp (- C_1 \lfloor k^{a + \varepsilon} \rfloor k^{- a})) + O (k^{a +
  \varepsilon} (k^{- 2 a + \varepsilon} + k^{- 2 \gamma + a +
  \varepsilon})) =  O (k^{- 2 \gamma + 2 a + \varepsilon} + k^{- a + \varepsilon}) ,
\end{align*}
where $O(\cdot)$ notation is used to ignore any constant terms to simplify the
notation. It immediately follows from \eqref{eq:vnkslkfjblkaf} that there exists
a constant $\mathscr{E}$ such that for any $1 - a > \varepsilon > 0$, the
following inequality holds:
\begin{align*}
  r_{k + 1} \leqslant & (k^{- 2 \gamma + 2 a + \varepsilon} + k^{- a +
  \varepsilon}) \mathscr{E}.
\end{align*}

\paragraph{Proof to \Cref{lemma:98620d3402422d993259e3a830cd49a398b44a20}}

First note that for any $\alpha_k \in [0, 1]$, the following
identity always holds for the difference between the estimation and the true
expectation. The first step comes from the update rule in \Cref{algo:theory}.
\begin{align*}
  & \tilde{e}_n^{(k + 1)} -\mathbb{E}_{\xi} \mathbf{e}_n (x^{(k)} ; \xi)\\
  = & (1 - \alpha_k)  \tilde{e}_n^{(k)} + \alpha_k \tilde{\mathbf{e}}_n
  (x^{(k)} ; \xi_k ; \hat{e}^{(k)}) -\mathbb{E}_{\xi} \mathbf{e}_n (x^{(k)} ;
  \xi)\\
  = & (1 - \alpha_k)  (\tilde{e}_n^{(k)} -\mathbb{E}_{\xi} \mathbf{e}_n
  (x^{(k)} ; \xi)) + \alpha_k  (\tilde{\mathbf{e}}_n (x^{(k)} ; \xi_k ;
  \hat{e}^{(k)}) -\mathbb{E}_{\xi} \mathbf{e}_n (x^{(k)} ; \xi))\\
  = & (1 - \alpha_k)  (\tilde{e}_n^{(k)} -\mathbb{E}_{\xi} \mathbf{e}_n (x^{(k
  - 1)} ; \xi)) + \alpha_k  (\tilde{\mathbf{e}}_n (x^{(k)} ; \xi_k ;
  \hat{e}^{(k)}) -\mathbb{E}_{\xi} \mathbf{e}_n (x^{(k)} ; \xi))\\
  & - (1 - \alpha_k)  (\mathbb{E}_{\xi} \mathbf{e}_n (x^{(k)} ; \xi)
  -\mathbb{E}_{\xi} \mathbf{e}_n (x^{(k - 1)} ; \xi)) . \numberthis
  \label{eq:vnzlksazajklsa}
\end{align*}
Then take $\ell_2$ norm on both sides. The second step comes from the fact
that $n$ is the last layer's index, so $\mathbf{e}_n (x ; \xi) = \hat{\mathbf{e}}_n (x^{(k)} ; \xi ;
  \hat{e}^{(k)})$.
\begin{align*}
  & \mathbb{E} \| \tilde{e}_n^{(k + 1)} -\mathbb{E}_{\xi} \mathbf{e}_n
  (x^{(k)} ; \xi)\|^2\\
  = & \mathbb{E} \left\| \begin{array}{c}
    (1 - \alpha_k) (\tilde{e}_n^{(k)} -\mathbb{E}_{\xi} \mathbf{e}_n (x^{(k -
    1)} ; \xi))\\
    + \alpha_k (\hat{\mathbf{e}}_n (x^{(k)} ; \xi_k ; \hat{e}^{(k)})
    -\mathbb{E}_{\xi} \mathbf{e}_n (x^{(k)} ; \xi))\\
    - (1 - \alpha_k) (\mathbb{E}_{\xi} \mathbf{e}_n (x^{(k)} ; \xi)
    -\mathbb{E}_{\xi} \mathbf{e}_n (x^{(k - 1)} ; \xi))
  \end{array} \right\|^2\\
  = & \mathbb{E} \left\| \begin{array}{c}
    (1 - \alpha_k) (\tilde{e}_n^{(k)} -\mathbb{E}_{\xi} \mathbf{e}_n (x^{(k -
    1)} ; \xi))\\
    + \alpha_k (\mathbf{e}_n (x^{(k)} ; \xi_k) -\mathbb{E}_{\xi} \mathbf{e}_n
    (x^{(k)} ; \xi))\\
    - (1 - \alpha_k) (\mathbb{E}_{\xi} \mathbf{e}_n (x^{(k)} ; \xi)
    -\mathbb{E}_{\xi} \mathbf{e}_n (x^{(k - 1)} ; \xi))
  \end{array} \right\|^2\\
  = & (1 - \alpha_k)^2 \mathbb{E} \left\| \begin{array}{c}
    \tilde{e}_n^{(k)} -\mathbb{E}_{\xi} \mathbf{e}_n (x^{(k - 1)} ; \xi)\\
    - (\mathbb{E}_{\xi} \mathbf{e}_n (x^{(k)} ; \xi) -\mathbb{E}_{\xi}
    \mathbf{e}_n (x^{(k - 1)} ; \xi))
  \end{array} \right\|^2 + \alpha_k^2 \mathbb{E} \| \mathbf{e}_n (x^{(k)} ;
  \xi_k) -\mathbb{E}_{\xi} \mathbf{e}_n (x^{(k)} ; \xi)\|^2\\
  & + 2 (1 - \alpha_k) \alpha_k \mathbb{E} \left\langle \begin{array}{c}
    \tilde{e}_n^{(k)} -\mathbb{E}_{\xi} \mathbf{e}_n (x^{(k - 1)} ; \xi)\\
    - (\mathbb{E}_{\xi} \mathbf{e}_n (x^{(k)} ; \xi) -\mathbb{E}_{\xi}
    \mathbf{e}_n (x^{(k - 1)} ; \xi))
  \end{array}, \mathbf{e}_n (x^{(k)} ; \xi_k) -\mathbb{E}_{\xi} \mathbf{e}_n
  (x^{(k)} ; \xi) \right\rangle\\
  = & (1 - \alpha_k)^2 \mathbb{E} \left\| \begin{array}{c}
    \tilde{e}_n^{(k)} -\mathbb{E}_{\xi} \mathbf{e}_n (x^{(k - 1)} ; \xi)\\
    - (\mathbb{E}_{\xi} \mathbf{e}_n (x^{(k)} ; \xi) -\mathbb{E}_{\xi}
    \mathbf{e}_n (x^{(k - 1)} ; \xi))
  \end{array} \right\|^2 + \alpha_k^2 \mathbb{E} \| \mathbf{e}_n (x^{(k)} ;
  \xi_k) -\mathbb{E}_{\xi} \mathbf{e}_n (x^{(k)} ; \xi)\|^2\\
  & + 2 (1 - \alpha_k) \alpha_k \mathbb{E} \left\langle \begin{array}{c}
    \tilde{e}_n^{(k)} -\mathbb{E}_{\xi} \mathbf{e}_n (x^{(k - 1)} ; \xi)\\
    - (\mathbb{E}_{\xi} \mathbf{e}_n (x^{(k)} ; \xi) -\mathbb{E}_{\xi}
    \mathbf{e}_n (x^{(k - 1)} ; \xi))
  \end{array}, \underbrace{\mathbb{E}_{\xi_k} \mathbf{e}_n (x^{(k)} ; \xi_k)
  -\mathbb{E}_{\xi} \mathbf{e}_n (x^{(k)} ; \xi)}_{= 0} \right\rangle\\
  \overset{\text{{\Cref{ass:a7f57b116c8ed2da7490c29a079fe3e51df929cc}}-2}}{\leqslant}
  & \underbrace{(1 - \alpha_k)^2 \mathbb{E} \left\| \begin{array}{c}
    \tilde{e}_n^{(k)} -\mathbb{E}_{\xi} \mathbf{e}_n (x^{(k - 1)} ; \xi)\\
    - (\mathbb{E}_{\xi} \mathbf{e}_n (x^{(k)} ; \xi) -\mathbb{E}_{\xi}
    \mathbf{e}_n (x^{(k - 1)} ; \xi))
  \end{array} \right\|^2}_{\backassign \mathscr{T}_1} + \alpha_k^2 \sigma^2 .
  \numberthis \label{eq:cvznklsajf}
\end{align*}
The second term is bounded by the bounded variance assumption. We denote the
first term as $\mathscr{T}_1$ and bound it separately.
With the fact that $\| x + y \|^2 \leqslant (1+\mathscr{d}) \|x\|^2 +
(1+\frac{1}{\mathscr{d}} ) \|y\|^2, \forall \mathscr{d}\in (0,1),\forall x,
\forall y$ we can bound $\mathscr{T}_1$ as following:
\begin{align*}
  \mathscr{T}_1 = & (1 - \alpha_k)^2 \mathbb{E} \left\| \begin{array}{c}
    (\tilde{e}_n^{(k)} -\mathbb{E}_{\xi} \mathbf{e}_n (x^{(k - 1)} ; \xi))\\
    - (\mathbb{E}_{\xi} \mathbf{e}_n (x^{(k)} ; \xi) -\mathbb{E}_{\xi}
    \mathbf{e}_n (x^{(k - 1)} ; \xi))
  \end{array} \right\|^2\\
  \leqslant & (1 - \alpha_k)^2  \left( (1 + \mathscr{d})\mathbb{E}\|
  \tilde{e}_n^{(k)} -\mathbb{E}_{\xi} \mathbf{e}_n (x^{(k - 1)} ; \xi)\|^2 +
  \left( 1 + \frac{1}{\mathscr{d}} \right) \mathbb{E}\|\mathbb{E}_{\xi}
  \mathbf{e}_n (x^{(k)} ; \xi) -\mathbb{E}_{\xi} \mathbf{e}_n (x^{(k - 1)} ;
  \xi)\|^2 \right), \forall \mathscr{d} \in (0, 1)\\
  \overset{\mathscr{d} \leftarrow \alpha_k}{=} & (1 - \alpha_k)^2  (1 +
  \alpha_k) \mathbb{E} \| \tilde{e}_n^{(k)} -\mathbb{E}_{\xi} \mathbf{e}_n
  (x^{(k - 1)} ; \xi)\|^2
   + \left( \frac{(1 + \alpha_k) (1 - \alpha_k)^2}{\alpha_k} \right)
  \mathbb{E} \|\mathbb{E}_{\xi} \mathbf{e}_n (x^{(k)} ; \xi) -\mathbb{E}_{\xi}
  e_n (x^{(k - 1)} ; \xi)\|^2\\
  \leqslant & (1 - \alpha_k) \mathbb{E} \| \tilde{e}_n^{(k)} -\mathbb{E}_{\xi}
  \mathbf{e}_n (x^{(k - 1)} ; \xi)\|^2 + \frac{1}{\alpha_k} \mathbb{E}
  \|\mathbb{E}_{\xi} \mathbf{e}_n (x^{(k)} ; \xi) -\mathbb{E}_{\xi}
  \mathbf{e}_n (x^{(k - 1)} ; \xi)\|^2\\
  \overset{\text{{\Cref{ass:a7f57b116c8ed2da7490c29a079fe3e51df929cc}}-4}}{\leqslant}
  & (1 - \alpha_k) \mathbb{E} \| \tilde{e}_n^{(k)} -\mathbb{E}_{\xi}
  \mathbf{e}_n (x^{(k - 1)} ; \xi)\|^2 + \frac{1}{\alpha_k} L\mathbb{E}
  \|x^{(k)} - x^{(k - 1)} \|^2 .\numberthis\label{eq:vjasdjmglajfldf}
\end{align*}
The last step comes from the Lipschitzian assumption on the functions.

Putting \eqref{eq:vjasdjmglajfldf} back into \eqref{eq:cvznklsajf} we obtain
\begin{align*}
  \mathbb{E} \| \tilde{e}_n^{(k + 1)} -\mathbb{E}_{\xi} \mathbf{e}_n (x^{(k)}
  ; \xi)\|^2 \leqslant & (1 - \alpha_k) \mathbb{E} \| \tilde{e}_n^{(k)}
  -\mathbb{E}_{\xi} \mathbf{e}_n (x^{(k - 1)} ; \xi)\|^2 + \frac{1}{\alpha_k}
  \underbrace{L\mathbb{E} \|x^{(k)} - x^{(k - 1)} \|^2}_{= O (\gamma_{k -
  1}^2)} + \alpha_k^2 \sigma^2, \numberthis
  \label{eq:7fecdb478ba95a9d31338d7d17aa8209fcd918a1}
\end{align*}
where the second term is of order $O (\gamma_{k - 1}^2)$ because the step length
is $\gamma_{k - 1}$ for the $(k - 1)$-th step and the gradient is bounded
according to {\Cref{ass:a7f57b116c8ed2da7490c29a079fe3e51df929cc}-1}. For the
case $i=n$, \eqref{eq:b5cf96af2f633135e356bd345b0fa5731e5e2b15} directly follows
from combining
\eqref{eq:7fecdb478ba95a9d31338d7d17aa8209fcd918a1},
\Cref{lemma:d083ec1ffe10dc6de84ee1d15cf12c068d8b28ab} and \Cref{ass:a7f57b116c8ed2da7490c29a079fe3e51df929cc}-6.

We then prove \eqref{eq:b5cf96af2f633135e356bd345b0fa5731e5e2b15} for all $i$ by
induction. Assuming for all $p> i$ and for any $1-a>\epsilon>0$ there exists a
constant $\mathscr{E}$ \eqref{eq:b5cf96af2f633135e356bd345b0fa5731e5e2b15} holds such that

\begin{equation}
    \mathbb{E} \| \tilde{e}_{p}^{(k + 1)} -\mathbb{E}_{\xi} \mathbf{e}_{p} (x^{(k)} ; \xi) \|^2 \leqslant \mathscr{E}(k^{- 2 \gamma + 2 a + \varepsilon} + k^{- a + \varepsilon}).\label{eq:hzzasdfmlk}
  \end{equation}

Similar to \eqref{eq:vnzlksazajklsa} we can split the difference between the estimation and the expectation into three parts:
\begin{align*}
  & \tilde{e}_i^{(k + 1)} -\mathbb{E}_{\xi} \mathbf{e}_i (x^{(k)} ; \xi)\\
  = & (1 - \alpha_k)  \tilde{e}_i^{(k)} + \alpha_k \hat{\mathbf{e}}_i
  (x^{(k)} ; \xi_k ; \hat{e}^{(k)}) -\mathbb{E}_{\xi} \mathbf{e}_i (x^{(k)} ;
  \xi)\\
  = & (1 - \alpha_k)  (\tilde{e}_i^{(k)} -\mathbb{E}_{\xi} \mathbf{e}_i (x^{(k
  - 1)} ; \xi)) + (1 - \alpha_k) \mathbb{E}_{\xi} \mathbf{e}_i (x^{(k - 1)} ;
  \xi) + \alpha_k \hat{\mathbf{e}}_i (x^{(k)} ; \xi_k ; \hat{e}^{(k)})
  -\mathbb{E}_{\xi} \mathbf{e}_i (x^{(k)} ; \xi)\\
  = & (1 - \alpha_k)  (\tilde{e}_i^{(k)} -\mathbb{E}_{\xi} \mathbf{e}_i (x^{(k
  - 1)} ; \xi)) + (1 - \alpha_k)  (\mathbb{E}_{\xi} \mathbf{e}_i (x^{(k - 1)}
  ; \xi) -\mathbb{E}_{\xi} \mathbf{e}_i (x^{(k)} ; \xi))\\
  & + \alpha_k  (\hat{\mathbf{e}}_i (x^{(k)} ; \xi_k ; \hat{e}^{(k)})
  -\mathbb{E}_{\xi} \mathbf{e}_i (x^{(k)} ; \xi)),
\end{align*}
Taking the $\ell_2$ norm on both sides we obtain
\begin{align*}
  & \mathbb{E} \| \tilde{e}_i^{(k + 1)} -\mathbb{E}_{\xi} \mathbf{e}_i
  (x^{(k)} ; \xi)\|^2\\
  \leqslant & \mathbb{E} \left\| \begin{array}{c}
    (1 - \alpha_k) (\tilde{e}_i^{(k)} -\mathbb{E}_{\xi} \mathbf{e}_i (x^{(k -
    1)} ; \xi))\\
    + (1 - \alpha_k) (\mathbb{E}_{\xi} \mathbf{e}_i (x^{(k - 1)} ; \xi)
    -\mathbb{E}_{\xi} \mathbf{e}_i (x^{(k)} ; \xi))\\
    + \alpha_k (\hat{\mathbf{e}}_i (x^{(k)} ; \xi_k ; \hat{e}^{(k)})
    -\mathbb{E}_{\xi} \mathbf{e}_i (x^{(k)} ; \xi))
  \end{array} \right\|^2\\
  = & \underbrace{(1 - \alpha_k)^2 \mathbb{E} \left\| \begin{array}{c}
    \tilde{e}_i^{(k)} -\mathbb{E}_{\xi} \mathbf{e}_i (x^{(k - 1)} ; \xi)\\
    + (\mathbb{E}_{\xi} \mathbf{e}_i (x^{(k - 1)} ; \xi) -\mathbb{E}_{\xi}
    \mathbf{e}_i (x^{(k)} ; \xi))
  \end{array} \right\|^2}_{\mathscr{T}_2}\\
  & + 2 (1 - \alpha_k) \alpha_k \underbrace{\left\langle \begin{array}{c}
    \tilde{e}_i^{(k)} -\mathbb{E}_{\xi} \mathbf{e}_i (x^{(k - 1)} ; \xi)\\
    +\mathbb{E}_{\xi} \mathbf{e}_i (x^{(k - 1)} ; \xi) -\mathbb{E}_{\xi}
    \mathbf{e}_i (x^{(k)} ; \xi)
  \end{array}, \hat{\mathbf{e}}_i (x^{(k)} ; \xi_k ; \hat{e}^{(k)})
  -\mathbb{E}_{\xi} \mathbf{e}_i (x^{(k)} ; \xi)
  \right\rangle}_{\mathscr{T}_4}\\
  & + \underbrace{\alpha_k^2 \mathbb{E} \| \mathbf{\hat{e}}_i (x^{(k)} ;
  \xi_k ; \hat{e}^{(k)}) -\mathbb{E}_{\xi} \mathbf{e}_i (x^{(k)} ;
  \xi)\|^2}_{\mathscr{T}_3},
\end{align*}
We define the there parts in the above inequality as $\mathscr{T}_2,
\mathscr{T}_3$ and $\mathscr{T}_4$ so that we can bound them one by one.
Firstly for $\mathscr{T}_3$, it can be easily bounded using \eqref{eq:hzzasdfmlk},
\Cref{ass:a7f57b116c8ed2da7490c29a079fe3e51df929cc}-2 and \Cref{ass:a7f57b116c8ed2da7490c29a079fe3e51df929cc}-4:
\begin{align*}
  \mathscr{T}_3 = & \alpha_k^2 \mathbb{E} \| \hat{\mathbf{e}}_i (x^{(k)} ;
  \xi_k ; \hat{e}^{(k)}) -\mathbb{E}_{\xi} \mathbf{e}_i (x^{(k)} ; \xi)\|^2\\
  \leqslant & \alpha_k^2 \mathbb{E} \| \hat{\mathbf{e}}_i (x^{(k)} ; \xi_k
  ; \hat{e}^{(k)}) -\mathbb{E}_{\xi} \hat{\mathbf{e}}_i (x^{(k)} ; \xi ;
  \hat{e}^{(k)})\|
   + \alpha_k^2 \mathbb{E} \|\mathbb{E}_{\xi} \hat{\mathbf{e}}_i
  (x^{(k)} ; \xi ; \hat{e}^{(k)}) -\mathbb{E}_{\xi} \mathbf{e}_i (x^{(k)} ;
  \xi)\|\\
  \overset{\text{{\Cref{ass:a7f57b116c8ed2da7490c29a079fe3e51df929cc}}-2}}{\leqslant} & \alpha_k^2 \sigma^2 + \alpha_k^2 \mathbb{E} \|\mathbb{E}_{\xi}
  \hat{\mathbf{e}}_i (x^{(k)} ; \xi ; \hat{e}^{(k)}) -\mathbb{E}_{\xi}
  \mathbf{e}_i (x^{(k)} ; \xi)\|\\
  \overset{\text{{\Cref{ass:a7f57b116c8ed2da7490c29a079fe3e51df929cc}}-4}}{\leqslant}
  & \alpha_k^2 \sigma^2 + \alpha_k^2 L \sum_{j = i+1}^{n} \mathbb{E} \|
  \tilde{e}_j^{(k)} -\mathbb{E}_{\xi} \mathbf{e}_j (x^{(k)} ; \xi)\|^2\\
  \overset{\eqref{eq:hzzasdfmlk}}{\leqslant} & \alpha_k^2 \sigma^2 + n \alpha_k^2 L \mathscr{E}.\numberthis \label{eq:zlammjqwqe}
\end{align*}
We then take a look at $\mathscr{T}_2$. It can be bounded similar to what we did
in \eqref{eq:vjasdjmglajfldf} to bound $\mathscr{T}_1$.
\begin{align*}
  \mathscr{T}_2 = & (1 - \alpha_k)^2 \mathbb{E} \left\| \begin{array}{c}
    \tilde{e}_i^{(k)} -\mathbb{E}_{\xi} \mathbf{e}_i (x^{(k - 1)} ; \xi)\\
    + (\mathbb{E}_{\xi} \mathbf{e}_i (x^{(k - 1)} ; \xi) -\mathbb{E}_{\xi}
    \mathbf{e}_i (x^{(k)} ; \xi))
  \end{array} \right\|^2\\
  \leqslant & (1 - \alpha_k)^2 \left( \begin{array}{c}
    (1 + \mathscr{d})\mathbb{E}\| \tilde{e}_i^{(k)} -\mathbb{E}_{\xi}
    \mathbf{e}_i (x^{(k - 1)} ; \xi)\|^2\\
    + \left( 1 + \frac{1}{\mathscr{d}} \right) \mathbb{E}\|\mathbb{E}_{\xi}
    \mathbf{e}_i (x^{(k - 1)} ; \xi) -\mathbb{E}_{\xi} \mathbf{e}_i (x^{(k)}
    ; \xi)\|^2
  \end{array} \right), \quad \forall \mathscr{d} \in (0, 1)\\
  \overset{\mathscr{d} \leftarrow \alpha_k}{\leqslant} & (1 - \alpha_k)
  \mathbb{E} \| \tilde{e}_i^{(k)} -\mathbb{E}_{\xi} \mathbf{e}_i (x^{(k - 1)}
  ; \xi)\|^2 + \frac{1}{\alpha_k} \mathbb{E} \|\mathbb{E}_{\xi} \mathbf{e}_i
  (x^{(k - 1)} ; \xi) -\mathbb{E}_{\xi} \mathbf{e}_i (x^{(k)} ; \xi)\|^2\\
  \overset{\text{{\Cref{ass:a7f57b116c8ed2da7490c29a079fe3e51df929cc}}-4}}{\leqslant}
  & (1 - \alpha_k) \mathbb{E} \| \tilde{e}_i^{(k)} -\mathbb{E}_{\xi}
  \mathbf{e}_i (x^{(k - 1)} ; \xi)\|^2 + \frac{L}{\alpha_k}
  (\mathbb{E}\|x^{(k)} - x^{(k - 1)} \|^2) .\numberthis\label{eq:zlaskfjgg}
\end{align*}
Finally we need to bound $\mathscr{T}_4$. This one is a litter harder than
$\mathscr{T}_2$ and $\mathscr{T}_3$. Different from what we were doing in
\eqref{eq:cvznklsajf}, we no longer have the nice property $\mathbf{e}_i (x ;
\xi) = \hat{\mathbf{e}}_i (x^{(k)} ; \xi ; \hat{e}^{(k)})$ since here we are
dealing with $i\neq n$.
We need to further split it and bound each part separately.
\begin{align*}
  \mathscr{T}_4 = & \left\langle \begin{array}{c}
    \tilde{e}_i^{(k)} -\mathbb{E}_{\xi} \mathbf{e}_i (x^{(k - 1)} ; \xi)\\
    +\mathbb{E}_{\xi} \mathbf{e}_i (x^{(k - 1)} ; \xi) -\mathbb{E}_{\xi}
    \mathbf{e}_i (x^{(k)} ; \xi)
  \end{array}, \mathbb{E}_{\xi} \hat{\mathbf{e}}_i (x^{(k)} ; \xi ;
  \hat{e}^{(k)}) -\mathbb{E}_{\xi} \mathbf{e}_i (x^{(k)} ; \xi)
  \right\rangle\\
  = & \underbrace{\langle \tilde{e}_i^{(k)} -\mathbb{E}_{\xi} \mathbf{e}_i
  (x^{(k - 1)} ; \xi), \mathbb{E}_{\xi} \hat{\mathbf{e}}_i (x^{(k)} ; \xi
  ; \hat{e}^{(k)}) -\mathbb{E}_{\xi} \mathbf{e}_i (x^{(k)} ; \xi)
  \rangle}_{\mathscr{T}_6}\\
  & + \underbrace{\langle \mathbb{E}_{\xi} \mathbf{e}_i (x^{(k - 1)} ; \xi)
  -\mathbb{E}_{\xi} \mathbf{e}_i (x^{(k)} ; \xi), \mathbb{E}_{\xi}
  \hat{\mathbf{e}}_i (x^{(k)} ; \xi ; \hat{e}^{(k)}) -\mathbb{E}_{\xi}
  \mathbf{e}_i (x^{(k)} ; \xi) \rangle}_{\mathscr{T}_7} . \numberthis
  \label{eq:cvjlzdkafml}
\end{align*}
We first bound $\mathscr{T}_7$:
\begin{align*}
  \mathscr{T}_7 = & \langle \mathbb{E}_{\xi} \mathbf{e}_i (x^{(k - 1)} ; \xi)
  -\mathbb{E}_{\xi} \mathbf{e}_i (x^{(k)} ; \xi), \mathbb{E}
  \widehat{\mathbf{e}}_i (x^{(k)} ; \xi_k ; \hat{e}^{(k)}) -\mathbb{E}_{\xi}
  \mathbf{e}_i (x^{(k)} ; \xi) \rangle\\
  \leqslant & \underbrace{\|\mathbb{E}_{\xi} \mathbf{e}_i (x^{(k - 1)} ; \xi)
  -\mathbb{E}_{\xi} \mathbf{e}_i (x^{(k)} ; \xi)\|}_{= O (\gamma_{k - 1})}
  \|\mathbb{E} \widehat{\mathbf{e}}_i (x^{(k)} ; \xi_k ; \hat{e}^{(k)})
  -\mathbb{E}_{\xi} \mathbf{e}_i (x^{(k)} ; \xi)\|\\
  = & O \left( \gamma_{k - 1} \|\mathbb{E} \widehat{\mathbf{e}}_i (x^{(k)} ;
  \xi_k ; \hat{e}^{(k)}) -\mathbb{E}_{\xi} \mathbf{e}_i (x^{(k)} ; \xi)\|
  \right)\numberthis \label{eq:valmermewqeq}\\
  \overset{\text{{\Cref{ass:a7f57b116c8ed2da7490c29a079fe3e51df929cc}}-4}}{\leqslant} & O \left( \gamma_{k - 1}  \sum_{j = i + 1}^n \| \hat{e}_j^{(k)}
  -\mathbb{E}_{\xi} \mathbf{e}_j (x^{(k)} ; \xi)\| \right)\\
  \overset{\eqref{eq:hzzasdfmlk}}{=} & O \left( \gamma_{k - 1}  \sqrt{k^{- 2 \gamma + 2 a + \varepsilon} + k^{-
  a + \varepsilon}} \right)
  = O (\gamma_{k - 1} (k^{- \gamma + a + \varepsilon / 2} + k^{- a / 2 +
  \varepsilon / 2}))\\
  = & O (k^{- 2 \gamma + a + \varepsilon / 2} + k^{- \gamma - a / 2 +
  \varepsilon / 2}), \numberthis\label{eq:zlkjsdkflagn}
\end{align*}
where in the second step we have $\|\mathbb{E}_{\xi} \mathbf{e}_i (x^{(k - 1)}
; \xi) -\mathbb{E}_{\xi} \mathbf{e}_i (x^{(k)} ; \xi)\| = O (\gamma_{k - 1})$
since by \Cref{ass:a7f57b116c8ed2da7490c29a079fe3e51df929cc}-4 we obtain $\|\mathbb{E}_{\xi} \mathbf{e}_i (x^{(k - 1)} ; \xi)
-\mathbb{E}_{\xi} \mathbf{e}_i (x^{(k)} ; \xi)\| \leqslant L \| x^{(k - 1)} -
x^{(k)} \|$. Then it follows from the same argument as in \eqref{eq:7fecdb478ba95a9d31338d7d17aa8209fcd918a1}.

After bounding $\mathscr{T}_7$, we start investigating $\mathscr{T}_6$. The last
step follows the same procedure as in \eqref{eq:valmermewqeq}.
\begin{align*}
  \mathscr{T}_6 = & \langle \tilde{e}_i^{(k)} -\mathbb{E}_{\xi} \mathbf{e}_i
  (x^{(k - 1)} ; \xi), \mathbb{E}\mathbf{e}_i (x^{(k)} ; \xi_k ;
  \hat{e}^{(k)}) -\mathbb{E}_{\xi} \mathbf{e}_i (x^{(k)} ; \xi) \rangle\\
  \leqslant & \frac{1}{2 \mathscr{d}_k}  \| \tilde{e}_i^{(k)}
  -\mathbb{E}_{\xi} \mathbf{e}_i (x^{(k - 1)} ; \xi)\|^2 + \mathscr{d}_k
  \|\mathbb{E}\mathbf{e}_i (x^{(k)} ; \xi_k ; \hat{e}^{(k)})
  -\mathbb{E}_{\xi} \mathbf{e}_i (x^{(k)} ; \xi)\|^2, \quad \forall
  \mathscr{d}_k > 0\\
  = & \frac{1}{2 \mathscr{d}_k}  \| \tilde{e}_i^{(k)} -\mathbb{E}_{\xi}
  \mathbf{e}_i (x^{(k - 1)} ; \xi)\|^2 + \mathscr{d}_k O (k^{- 2 \gamma + 2 a
  + \varepsilon} + k^{- a + \varepsilon}),\quad \forall \mathscr{d}_k > 0.\numberthis \label{eq:zpppsafjasdf}
\end{align*}

Finally plug \eqref{eq:zlammjqwqe}, \eqref{eq:zlaskfjgg},
  \eqref{eq:zpppsafjasdf} and \eqref{eq:zlkjsdkflagn} back into \eqref{eq:cvjlzdkafml}. By choosing
  $\mathscr{d}_k$ in \eqref{eq:zpppsafjasdf} to be 1 we obtain:
\begin{align*}
  & \mathbb{E} \| \tilde{e}_i^{(k + 1)} -\mathbb{E}_{\xi} \mathbf{e}_i
  (x^{(k)} ; \xi)\|^2\\
  \leqslant & \alpha_k^2 \sigma^2 + n \alpha_k^2 L \mathscr{E}\\
  & + (1 - \alpha_k) \mathbb{E} \| \tilde{e}_i^{(k)} -\mathbb{E}_{\xi}
  \mathbf{e}_i (x^{(k - 1)} ; \xi)\|^2 + \frac{L}{\alpha_k}
  (\mathbb{E}\|x^{(k)} - x^{(k - 1)} \|^2)\\
  & + 2 (1 - \alpha_k) \alpha_k O (k^{- 2 \gamma + a + \varepsilon / 2} +
  k^{- \gamma - a / 2 + \varepsilon / 2})\\
  & + (1 - \alpha_k) \alpha_k  \frac{1}{\mathscr{d}_k}  \| \tilde{e}_i^{(k)}
  -\mathbb{E}_{\xi} \mathbf{e}_i (x^{(k - 1)} ; \xi)\|^2 + 2 (1 - \alpha_k)
  \alpha_k  \mathscr{d}_k O (k^{- 2 \gamma + 2 a + \varepsilon} + k^{- a +
  \varepsilon})\\
  \overset{\mathscr{d}_k \leftarrow 1}{\leqslant} & \alpha_k^2 \sigma^2 + n
  \alpha_k^2 L \mathscr{E}\\
  & + (1 - \alpha_k) \mathbb{E} \| \tilde{e}_i^{(k)} -\mathbb{E}_{\xi}
  \mathbf{e}_i (x^{(k - 1)} ; \xi)\|^2 + \frac{L}{\alpha_k}
  (\mathbb{E}\|x^{(k)} - x^{(k - 1)} \|^2)\\
  & + 2 (1 - \alpha_k) \alpha_k O (k^{- 2 \gamma + a + \varepsilon / 2} +
  k^{- \gamma - a / 2 + \varepsilon / 2})\\
  & + (1 - \alpha_k) \alpha_k  \frac{1}{2}  \| \tilde{e}_i^{(k)}
  -\mathbb{E}_{\xi} \mathbf{e}_i (x^{(k - 1)} ; \xi)\|^2 + 4 (1 - \alpha_k)
  \alpha_k O (k^{- 2 \gamma + 2 a + \varepsilon} + k^{- a + \varepsilon})\\
  \leqslant & \left( 1 - \frac{\alpha_k}{2} \right) \mathbb{E} \|
  \tilde{e}_i^{(k)} -\mathbb{E}_{\xi} \mathbf{e}_i (x^{(k - 1)} ; \xi)\|^2 +
  \alpha_k^2 \sigma^2 + n \alpha_k^2 L \mathscr{E} + \frac{L \gamma_{k -
  1}^2}{\alpha_k}  \mathscr{G}^2\\
  & + 2 (1 - \alpha_k) \alpha_k O (k^{- 2 \gamma + a + \varepsilon / 2} +
  k^{- \gamma - a / 2 + \varepsilon / 2}) + 4 (1 - \alpha_k) \alpha_k O (k^{-
  2 \gamma + 2 a + \varepsilon} + k^{- a + \varepsilon})\\
  = & \left( 1 - \frac{\alpha_k}{2} \right) \mathbb{E} \| \tilde{e}_i^{(k)}
  -\mathbb{E}_{\xi} \mathbf{e}_i (x^{(k - 1)} ; \xi)\|^2 + \alpha_k^2
  \sigma^2\\
  & + n \alpha_k^2 L \mathscr{E} + \frac{L \gamma_{k - 1}^2}{\alpha_k}
  \mathscr{G}^2 + O (k^{- 2 \gamma + a + \varepsilon} + k^{- 2 a +
  \varepsilon}) . \numberthis
  \label{eq:fd18e7bf2a363150a68972311612531f1b47f604}
\end{align*}

By combining \eqref{eq:fd18e7bf2a363150a68972311612531f1b47f604} and \Cref{lemma:d083ec1ffe10dc6de84ee1d15cf12c068d8b28ab} we obtain \eqref{eq:b5cf96af2f633135e356bd345b0fa5731e5e2b15}.
\eqref{eq:2ab3f273dad20de3dbd553953676c02229427496} follows from
\eqref{eq:7fecdb478ba95a9d31338d7d17aa8209fcd918a1} and
\eqref{eq:fd18e7bf2a363150a68972311612531f1b47f604}.

\paragraph{Proof to \Cref{thm:384575d2d125c25265b5d7bf9c978acb95307110}}
It directly follows from the definition of Lipschitz condition of $f$ that,
\begin{align*}
  f (x^{(k + 1)}) \overset{\text{\Cref{ass:a7f57b116c8ed2da7490c29a079fe3e51df929cc}-4}}{\leqslant} & f
  (x^{(k)}) + \langle \partial f (x^{(k)}), x^{(k + 1)} - x^{(k)} \rangle + L
  \| x^{(k + 1)} - x^{(k)} \|^2\\
  = & f (x^{(k)}) + \langle \partial f (x^{(k)}), - \gamma_k g^{(k)} \rangle +
  L \| x^{(k + 1)} - x^{(k)} \|^2\\
  = & f (x^{(k)}) - \gamma_k  \| \partial f (x^{(k)}) \|^2
   + \underbrace{\langle \partial f (x^{(k)}), - \gamma_k  (g^{(k)} -
  \partial f (x^{(k)})) \rangle}_{\mathscr{T}\text{cross}} + L \underbrace{\|
  x^{(k + 1)} - x^{(k)} \|^2}_{\mathscr{T}\text{progress}}.\numberthis \label{eq:cnzldfsalnvasa}
\end{align*}
Here we define two new terms and try to bound the expectation of
$\mathscr{T}\text{cross}$ and $\mathscr{T}\text{progress}$ separately as shown
below:
\begin{align*}
  \mathbb{E} \mathscr{T} \text{cross} = & \mathbb{E} \langle \partial f
  (x^{(k)}), - \gamma_k (g^{(k)} - \partial f (x^{(k)})) \rangle\\
  \leqslant & \frac{2 \gamma_k^2 L_g}{\alpha_{k + 1}} \mathbb{E} \| \partial f
  (x^{(k)})\|^2 + \frac{\alpha_{k + 1}}{2 L_g} \mathbb{E} \|g^{(k)} - \partial
  f (x^{(k)})\|^2\\
  \overset{\text{{\Cref{ass:a7f57b116c8ed2da7490c29a079fe3e51df929cc}}-3}}{\leqslant}
  & \frac{2 \gamma_k^2 L_g}{\alpha_{k + 1}} \mathbb{E} \| \partial f
  (x^{(k)})\|^2 + \frac{\alpha_{k + 1}}{2}  \sum_{i = 1}^n \| \tilde{e}_i^{(k
  + 1)} -\mathbb{E}_{\xi} [\mathbf{e}_i (x^{(k)} ; \xi)]\|^2, \numberthis
  \label{eq:czlnsdklnqw}
\end{align*}
and
\begin{align*}
  \mathbb{E}\mathscr{T}\text{progress}= & \mathbb{E} \| x_{k + 1} - x_k
  \|^2
  \overset{\text{\Cref{ass:a7f57b116c8ed2da7490c29a079fe3e51df929cc}-1}}{\leqslant} \gamma^2_k \mathscr{G}^2 . \numberthis \label{eq:vnklarjw}
\end{align*}

With the help of \eqref{eq:czlnsdklnqw} and \eqref{eq:vnklarjw}, \eqref{eq:cnzldfsalnvasa} becomes
\begin{align*}
  \mathbb{E}f (x^{(k + 1)}) \leqslant & \mathbb{E}f (x^{(k)}) - \gamma_k
  \mathbb{E} \| \partial f (x^{(k)}) \|^2\\
  & + \frac{2 \gamma_k^2 L_g}{\alpha_{k + 1}} \mathbb{E} \| \partial f
  (x^{(k)}) \|^2 + \frac{\alpha_{k + 1}}{2}  \sum_{i = 1}^{n - 1} \mathbb{E}
  \| \tilde{e}_i^{(k + 1)} -\mathbb{E}_{\xi_k} [e_i (x^{(k)} ; \xi_k)] \|^2 +
  \gamma_k^2 L\mathscr{G}^2 . \numberthis \label{eq:nzlfsdkafa}
\end{align*}

To show how this converges, we define the following $\mathscr{d}_k $ to derive a
recursive relation for \eqref{eq:nzlfsdkafa}:
\begin{align*}
  \mathscr{d}_k := & f (x^{(k)}) + \sum_{i = 1}^{n - 1} \| \tilde{e}_i^{(k +
  1)} -\mathbb{E}_{\xi_k} [e_i (x^{(k)} ; \xi_k)] \|^2.
\end{align*}

Then the recursive relation is derived by observing
\begin{align*}
  \mathbb{E}\mathscr{d}_{k + 1} \leqslant & \mathbb{E}f (x^{(k)}) - \gamma_k
  \mathbb{E} \| \partial f (x^{(k)}) \|^2
   + \frac{2 \gamma_k^2 L_g}{\alpha_{k + 1}} \mathbb{E} \| \partial f
  (x^{(k)}) \|^2 + \frac{\alpha_{k + 1}}{2}  \sum_{i = 1}^{n - 1} \|
  \tilde{e}_i^{(k + 1)} -\mathbb{E}_{\xi_k} [e_i (x^{(k)} ; \xi_k)] \|^2 +
  \gamma_k^2 L\mathscr{G}^2\\
  & + \sum_{i = 1}^{n - 1} \| \tilde{e}_i^{(k + 2)} -\mathbb{E}_{\xi_{k + 1}}
  [e_i (x^{(k + 1)} ; \xi_{k + 1})] \|^2\\
  \overset{\text{\Cref{lemma:98620d3402422d993259e3a830cd49a398b44a20}}}{\leqslant} & \mathbb{E}f (x^{(k)}) - \gamma_k
  \mathbb{E} \| \partial f (x^{(k)}) \|^2\\
  & + \frac{2 \gamma_k^2 L_g}{\alpha_{k + 1}} \mathbb{E} \| \partial f
  (x^{(k)}) \|^2 + \frac{\alpha_{k + 1}}{2}  \sum_{i = 1}^{n - 1} \mathbb{E}
  \| \tilde{e}_i^{(k + 1)} -\mathbb{E}_{\xi_k} [e_i (x^{(k)} ; \xi_k)] \|^2 +
  \gamma_k^2 L\mathscr{G}^2\\
  & + \sum_{i = 1}^{n - 1} \left( \left( 1 - \frac{\alpha_{k + 1}}{2} \right)
  \mathbb{E} \|  \tilde{e}_i^{(k + 1)} -\mathbb{E}_{\xi} \mathbf{e}_i (x^{(k + 1)} ;
  \xi) \|^2 +\mathscr{C} ((k + 1)^{- 2 \gamma + a + \varepsilon} + (k + 1)^{-
  2 a + \varepsilon}) \right)\\
  \leqslant & \mathbb{E}f (x^{(k)}) - \gamma_k \mathbb{E} \| \partial f
  (x^{(k)}) \|^2
   + \frac{2 \gamma_k^2 L_g}{\alpha_{k + 1}} \mathbb{E} \| \partial f
  (x^{(k)}) \|^2 + \sum_{i = 1}^{n - 1} \mathbb{E} \|  \tilde{e}_i^{(k + 1)}
  -\mathbb{E}_{\xi} \mathbf{e}_i (x^{(k)} ; \xi) \|^2
   + O (k^{- 2 \gamma + a + \varepsilon} + k^{- 2 a + \varepsilon})\\
  = & \mathbb{E}\mathscr{d}_k - \left( \gamma_k - \frac{2 \gamma_k^2
  L_g}{\alpha_{k + 1}} \right) \mathbb{E} \| \partial f (x^{(k)}) \|^2 + O
  (k^{- 2 \gamma + a + \varepsilon} + k^{- 2 a + \varepsilon}) .
\end{align*}

Thus as long as $a-2\gamma < -1, a < 1/2$ we can choose $\epsilon$ such that there
exists a constant $\mathscr{R}$:
\begin{align*}
  \mathbb{E}\mathscr{d}_{K + 1} \leqslant & \mathbb{E}\mathscr{d}_0 -
  \sum_{k = 0}^K \left( \gamma_k - \frac{2 \gamma_k^2 L_g}{\alpha_{k + 1}}
  \right) \mathbb{E} \| \partial f (x^{(k)}) \|^2 +\mathscr{R}\\
  \sum_{k = 0}^K \left( \gamma_k - \frac{2 \gamma_k^2
  L_g}{\alpha_{k + 1}} \right) \mathbb{E} \| \partial f (x^{(k)}) \|^2
  \leqslant & \mathscr{R}+\mathbb{E}\mathscr{d}_0 -\mathbb{E}\mathscr{d}_{K
  + 1}\\
  \overset{\frac{\gamma_k L_g}{\alpha_{k+1}} \leqslant \frac{1}{2}}{\longrightarrow} \frac{\sum_{k =
  0}^K \gamma_k \mathbb{E} \| \partial f (x^{(k)}) \|^2}{\sum_{k = 0}^K
  \gamma_k} \leqslant & \frac{2 (\mathscr{R}+\mathbb{E}\mathscr{d}_0
  -\mathbb{E}\mathscr{d}_{K + 1})}{\sum_{k = 0}^K \gamma_k} .
\end{align*}
By \Cref{ass:a7f57b116c8ed2da7490c29a079fe3e51df929cc}-5 we complete the proof.

\paragraph{Proof to \Cref{coro:8dce0d0cfc00f8fa3138e18c07037ccfc2ffd86a}}

With the given choice of $\alpha_k$ and $\gamma_k$ the prerequisites in
\Cref{thm:384575d2d125c25265b5d7bf9c978acb95307110} are satisfied.

Then by \Cref{thm:384575d2d125c25265b5d7bf9c978acb95307110} there exists a
constant $\mathscr{H}$ (which may differ from the constant in
\Cref{thm:384575d2d125c25265b5d7bf9c978acb95307110}) such that
\begin{align*}
  \sum_{k = 0}^K (k + 2)^{- 4 / 5} \mathbb{E} \| \partial f (x^{(k)}) \|^2
  \leqslant & \mathscr{H}\\
  \sum_{k = 0}^K \mathbb{E} \| \partial f (x^{(k)}) \|^2
  \leqslant & \frac{\mathscr{H}}{(K + 2)^{- 4 / 5}}\\
  \frac{\sum_{k = 0}^K \mathbb{E} \| \partial f (x^{(k)})
  \|^2}{K + 2} \leqslant & \frac{\mathscr{H}}{(K + 2)^{1 / 5}} ,
\end{align*}
completing the proof.

\end{document}